\documentclass[a4paper,12pt,reqno]{amsart}
\usepackage{vmargin}
\usepackage[pdftex]{graphicx}
\usepackage{amsfonts}
\usepackage{amssymb}
\usepackage{mathrsfs}
\usepackage{stmaryrd}
\usepackage{amsmath}
\usepackage{amsthm}
\usepackage[shortlabels]{enumitem}
\usepackage{nomencl}
\usepackage[bookmarks=true]{hyperref}   
\usepackage{esint}
\usepackage{dsfont}
\usepackage{upgreek}
\usepackage{xcolor}
\usepackage{slashed}
\usepackage{scalerel}
\usepackage{comment}
\usepackage[utf8]{inputenc}

\renewcommand{\imath}{{\rm i}}	

\makeatletter
\renewcommand\section{\@startsection{section}{1}{\z@}%
                       {-3\p@ \@plus -4\p@ \@minus -4\p@}%
                       {3\p@ \@plus 4\p@ \@minus 4\p@}%
                      {\normalfont\normalsize\centering\scshape}}
\makeatother

\newcommand{\mynote}[1]{#1}
\renewcommand{\mynote}[1]{}		

\numberwithin{equation}{section} 	

\hypersetup{
    colorlinks,%
}

\newcommand{\Rosen}{Ros\'en}

\author{Lashi Bandara}
\author{Andreas Rosén}
\title[Riesz continuity under perturbations of local boundary conditions]
{Riesz continuity of the Atiyah-Singer Dirac operator under perturbations of local boundary conditions}
\date{\today}

\address{Lashi Bandara, 
Institut für Mathematik,
Universität Potsdam, 
D-14476, Potsdam OT Golm, Germany
}
\urladdr{\href{http://www.math.uni-potsdam.de/~bandara}{http://www.math.uni-potsdam.de/~bandara}}
\email{\href{mailto:lashi.bandara@uni-potsdam.de}{lashi.bandara@uni-potsdam.de}}

\address{Andreas \Rosen, Mathematical Sciences,
Chalmers University of Technology and University of Gothenburg, SE-412 96, Gothenburg, Sweden}
\urladdr{\href{http://www.math.chalmers.se/~rosenan}{http://www.math.chalmers.se/~rosenan}}
\email{\href{mailto:andreas.rosen@chalmers.se}{andreas.rosen@chalmers.se}}

\dedicatory{Dedicated to the memory of Alan G. R. McIntosh}

\keywords{Riesz continuity, Dirac operator, spectral flow, 
functional calculus, boundary value problems, real-variable harmonic analysis}
\subjclass[2010]{58J05, 58J32, 58J37, 58J30, 42B37, 35J46, 35J56}

\def\colour{\colour}

\def\colour{\color}

\newtheorem{theorem}{Theorem}[section]

\newtheorem{lemma}[theorem]{Lemma}
\newtheorem{proposition}[theorem]{Proposition}

\newtheorem{remark}[theorem]{Remark}

\newtheorem{example}[theorem]{Example}

\newcommand{\mdot}{\cdotp}

\newcommand{\cbrac}[1]{\left(#1\right)}

\newcommand{\dbrac}[1]{\left\{#1\right\}}
\newcommand{\modulus}[1]{|#1|}
\newcommand{\lmodulus}[1]{\left|#1\right|}
\newcommand{\set}[1]{\dbrac{#1}}

\newcommand{\dom}{ {\mathcal{D}}}

\newcommand{\comp}{\, \circ\, }

\newcommand{\e}{\mathrm{e}}

\newcommand{\R}{\mathbb{R}}
\newcommand{\C}{\mathbb{C}}

\newcommand{\Na}{\ensuremath{\mathbb{N}}}

\newcommand{\script}[1]{\mathscr{#1}}

\DeclareMathOperator{\re}{Re}			

\renewcommand{\emptyset}{\varnothing}
\newcommand{\union}{\cup}
\newcommand{\Union}{\bigcup}
\newcommand{\intersect}{\cap}

\newcommand{\disunion}{\sqcup}

\newcommand{\rest}[1]{{{\lvert_{}}_{}}_{#1}}
\newcommand{\close}[1]{\overline{#1}}		

\newcommand{\ind}[1]{\raisebox{\depth}{\(\chi\)}_{#1}}	

\renewcommand{\epsilon}{\varepsilon}

\renewcommand{\phi}{\varphi}

\newcommand{\ch}[1]{\upchi_{{#1}}}	

\newcommand{\embed}{\hookrightarrow}		

\newcommand{\tensor}{\otimes}

\newcommand{\comm}[1]{[#1]}		

\newcommand{\End}{\mathrm{End}}

\newcommand{\norm}[1]{\| #1 \|}			
\newcommand{\snorm}[1]{\left\| #1 \right\|}			

\newcommand{\spt}[1]{{\rm spt} {\text{ }}#1}	

\DeclareMathOperator{\gap}{\hat{\updelta}}		
\newcommand{\cinterpol}[2]{{[#1]}_{#2}}			

\newcommand{\conform}{{\upomega}}		

\DeclareMathOperator{\tr}{tr}			

\DeclareMathOperator{\len}{\ell}			
\DeclareMathOperator{\rad}{rad}			

\DeclareMathOperator{\divv}{div}		

\newcommand{\Ric}{{\rm Ric}}			

\DeclareMathOperator{\inj}{inj} 		

\newcommand{\bnd}{\partial}			
\newcommand{\bdy}{\partial}			
\newcommand{\interior}[1]{\mathring{#1}}	
\newcommand{\intr}[1]{\interior{#1}}			
\newcommand{\Forms}[1][{}]{\mathbf{\Omega}^{#1}}		
\newcommand{\Tensors}[1][{}]{{\mathcal{T}}^{(#1)}}	
\newcommand{\Sect}{\mathbf{\Gamma}}		
\newcommand{\tanb}{{\rm T}}		
\newcommand{\cotanb}{{\rm T}^\ast}	

\newcommand{\pullb}[1]{#1^\ast}			

\newcommand{\on}{\vec{\mathrm{n}}}

\DeclareFontFamily{OT1}{restrictfont}{}
\DeclareFontShape{OT1}{restrictfont}{m}{n}{<-> fmvr8x}{}

\newcommand{\adj}[1]{{#1}^\ast}			
\newcommand{\extd}{{\rm d}}			

\newcommand{\inprod}[1]{\left\langle #1 \right\rangle}	
\newcommand{\conn}[1][{}]{{\nabla_{{#1}}}}		
\newcommand{\Leb}[1][{}]{\script{L}^{#1}}			

\newcommand{\Cliff}[1][{}]{\Delta^{#1}}		
\DeclareMathOperator{\Spin}{Spin}			
\DeclareMathOperator{\Spinors}{\slashed{\Delta}}	
\newcommand{\spin}[1]{\slashed{#1}}		
\newcommand{\Prin}[1]{\mathrm{P}_{#1}}		

\newcommand{\mi}{\mathrm{min}} 
\newcommand{\mx}{\mathrm{max}}
\newcommand{\sym}{\mathrm{sym}}
\newcommand{\bddlf}{\mathcal{L}} 	

\newcommand{\Lp}[2][{}]{{\rm L}^{#2}_{\rm #1}}		
\newcommand{\Ck}[2][{}]{{\rm C}^{#2}_{\rm #1}}		
\newcommand{\Hard}[2][{}]{{\rm H}^{#2}_{\rm #1}}		
\newcommand{\SobH}[2][{}]{\Hard[#1]{\rm #2}}
\newcommand{\Lips}[1][{}]{{\rm Lip}_{\rm #1}}		

\DeclareMathOperator{\Tr}{\mathscr{R}}				

\newcommand{\Sec}[1]{\mathrm{S}_{#1}}
\newcommand{\OSec}[1]{\mathrm{S}^\mathrm{o}_{#1}}

\newcommand{\maxx}[1]{\left<#1\right>}

\newcommand{\iden}{{\mathrm{I}}}

\newcommand{\Lap}{\Delta}			

\newcommand{\Q}[1][{}]{Q_{#1}}			
\newcommand{\DyQ}{\script{Q}}			
\newcommand{\ancester}[1]{\widehat{#1}}		

\newcommand{\brad}{\rho} 			
\newcommand{\scale}{\mathrm{t_S}}			
\newcommand{\jscale}{\mathrm{J}}

\newcommand{\sC}{\script{C}}

\newcommand{\sW}{\script{W}}

\newcommand{\cA}{\mathcal{A}}
\newcommand{\cB}{\mathcal{B}}

\newcommand{\cE}{\mathcal{E}}
\newcommand{\cV}{\mathcal{V}}
\newcommand{\cM}{\mathcal{M}} 

\newcommand{\cP}{\mathcal{P}}
\newcommand{\cW}{\mathcal{W}}

\newcommand{\Hol}{{\rm Hol}}

\newcommand{\mg}{\mathrm{g}}
\newcommand{\mgt}{\tilde{\mg}}
\newcommand{\mh}{\mathrm{h}}

\newcommand{\Poincare}{Poincar\'e~}		

\newcommand{\Av}{\mathbb{E}}
\newcommand{\Pri}{\upgamma}

\newcommand{\CBox}{\mathrm{R}}

\newcommand{\B}{\mathrm{B}}
\newcommand{\A}{\mathrm{A}}
\newcommand{\U}{\mathrm{U}}

\newcommand{\Dir}{{\rm D} }
\newcommand{\Dirb}{\Dir}		
\newcommand{\Dirp}{\tilde{\Dir}}

\newcommand{\BDir}{\slashed{\partial}}

\newcommand{\Qqb}{{\rm Q}}
\newcommand{\Ppb}{{\rm P}}
\newcommand{\Qq}{\tilde{\Qqb}}
\newcommand{\Pp}{\tilde{\Ppb}}
\newcommand{\QQ}{{\rm\bf Q}}

\newcommand{\Rrb}{{\rm R}}
\newcommand{\Rr}{\tilde{\Rrb}}
\newcommand{\Uub}{{\rm U}}
\newcommand{\Uu}{\tilde{\Uub}}

\newcommand{\dtt}[1][{t}]{\frac{d#1}{#1}}

\newcommand{\const}{\mathrm{C}} 
\newcommand{\met}{\uprho}		

\newcommand{\Ball}{\mathrm{B}}

\newcommand{\SDir}{\slashed{\Dir}}

\newcommand{\Sconn}[1][{}]{{{\nabla}_{{#1}}}}		

\newcommand{\rep}{\cdot}

\newcommand{\RNum}[1]{\uppercase\expandafter{\romannumeral #1\relax}}

\begin{document}

\maketitle
\vspace*{-2em}
\begin{abstract}
On a smooth complete Riemannian spin manifold with smooth
compact boundary,
we demonstrate that the Atiyah-Singer Dirac operator
$\SDir_\cB$ in $\Lp{2}$ depends Riesz continuously
on $\Lp{\infty}$ perturbations of local boundary 
conditions $\cB$.
The Lipschitz bound for the map $\cB \to \SDir_\cB(1 + \SDir_\cB^2)^{-\frac{1}{2}}$
depends on Lipschitz smoothness and ellipticity of $\cB$ and bounds on
Ricci curvature and its first derivatives
as well as a lower bound on injectivity radius  away from a compact neighbourhood of the boundary.
 More generally, we prove perturbation estimates 
for functional calculi of elliptic operators 
on manifolds with local boundary conditions.
\end{abstract}
\tableofcontents
\vspace*{-2em}

\parindent0cm
\setlength{\parskip}{\baselineskip}

\section{Introduction}

The aim of this paper and its companion \cite{BMcR} has been to prove perturbation estimates of quantities of the form
$$
  \snorm{\frac{\Dirp}{\sqrt{\iden+\Dirp^2}} - \frac{\Dirb}{\sqrt{\iden+\Dirb^2}}}_{\Lp{2}(\cM,\cV)\to \Lp{2}(\cM,\cV)},
$$ 
where $\Dirb$ and $\Dirp$ are self-adjoint elliptic first-order partial differential operators, acting 
on sections of a vector bundle $\cV$ over a smooth manifold $\cM$.
The symbol $f(\zeta)=\zeta(1+\zeta^2)^{-\frac{1}{2}}$ is a motivating example, yielding continuity results in the Riesz sense, but 
our methods apply equally well to more general holomorphic symbols around $\R$, 
which may be discontinuous at $\infty$.
In \cite{BMcR}, together with Alan McIntosh, we obtained results on complete manifolds  $(\cM,\mg)$ 
without boundary.
In that case, the main example of operators $\Dirb$ and $\Dirp$ was the Atiyah--Singer Dirac operators on $\cM$ with respect to two 
different metrics $\mg$ and $\mgt$.
The bound obtained was
$$
\snorm{\frac{\Dirp}{\sqrt{\iden+\Dirp^2}} - \frac{\Dirb}{\sqrt{\iden+\Dirb^2}}}_{\Lp{2}(\cM,\cV)\to \Lp{2}(\cM,\cV)}
 \lesssim \norm{\mgt -\mg}_{\Lp{\infty}(\Tensors[2,0]\cM)}, 
$$ 
where the implicit constant depends on certain geometric quantities. Note that the two Dirac operators themselves depend also on the first derivatives of the metrics.

In the present paper, we consider the corresponding perturbation estimate on a manifold $\cM$ (possibly noncompact) 
with smooth, compact boundary $\Sigma= \bdy\cM$.
Our motivating example in this case is when both $\Dirb$ and $\Dirp$ are the Atiyah--Singer Dirac operator, 
but with respect to two different local boundary conditions, defined through two different subbundles $\cE$ and $\tilde\cE$ of $\cV|_\Sigma$. 
For each boundary condition we assume self-adjointness and ellipticity so that the domains of $\Dirb$ and $\Dirp$ are closed subspaces of $\SobH{1}(\cV)$.
The bound we obtain is
\begin{equation}  \label{eq:mainpertest}
\snorm{\frac{\Dirp}{\sqrt{\iden+\Dirp^2}} - \frac{\Dirb}{\sqrt{\iden+\Dirb^2}}}_{\Lp{2}(\cM,\cV)\to \Lp{2}(\cM,\cV)} 
	\lesssim \norm{\gap(\tilde{\cE}_x, \cE_x )}_{\Lp{\infty}(\Sigma)},
\end{equation}
 where $\gap(\cE_x, \tilde{\cE}_x) = \modulus{\pi_{\cE}(x) - \pi_{\tilde{\cE}}(x)}$
and $\pi_{\cE}$ and $\pi_{\tilde{\cE}}$ are the orthogonal projectors
from $\cV\rest{\Sigma}$ to $\cE$ and $\tilde{\cE}$ respectively.  
Again the implicit constant in the estimate depends on a number of geometric quantities which 
 we list completely.

As described in the introduction of \cite{BMcR}, an important application of these perturbation
estimates is the study of spectral flow for unbounded self-adjoint operators. 
The study of the spectral flow was initiated by Atiyah and Singer in \cite{AS69} and 
has important connections
to particle physics. An analytic formulation of the spectral flow was given by Phillips in \cite{P96}
and typically, the \emph{gap metric} 
$$ \snorm{ \frac{\imath + \Dirp}{\imath - \Dirp} - \frac{\imath + \Dirb}{\imath - \Dirb}}_{\Lp{2}(\cM; \cV) \to  \Lp{2}(\cM;\cV)}$$
is used to understand the spectral flow for unbounded operators.
The Riesz topology is a preferred alternative
since  the spectral flow in this topology better connects to topological 
and $K$-theoretic aspects of the spectral flow, which were
observed in \cite{AS69} for the case of bounded
self-adjoint Fredholm operators.
The main disadvantage is that it is typically harder to establish 
continuity in the Riesz topology. 
In particular we refer to the open problem pointed out by Lesch in the introduction of \cite{Lesch}, 
namely whether a Dirac operator on a compact manifold with boundary depends Riesz continuously on 
pseudo-differential boundary conditions imposed on the operator.
 
The present paper answers this questions to the positive, in the special case of local boundary conditions. 
Self-adjoint local boundary conditions are typically physical 
and a very large subclass of the so-called \emph{Chiral} conditions
are listed in \cite{HMR}  by Hijazi, Montiel and Roldán as being self-adjoint boundary conditions. 
In particular, these exist in even dimensions or when the manifold is a space-like hypersurface in spacetime.   
The case of non-local boundary conditions defined by pseudo-differential projections appear to be beyond the scope of the methods used in the present paper but we anticipate they will be the object of further investigations in the future.
The local nature of the boundary conditions enter the proof in a number of instances, but the most serious occurrence concern the so-called exponential off-diagonal estimates, 
which relies on the domains of the operators being preserved under multiplication by smooth, bounded functions. 
It is important to note that the right hand sides in the perturbation estimates that we obtain, namely 
$\norm{\mgt -\mg}_{\Lp{\infty}(\Tensors[2,0]\cM)}$ and $\norm{\gap(\tilde{\cE}_x, \cE_x)}_{\Lp{\infty}(\Sigma)}$, 
are supremum norms,  which are smaller than estimates that can be obtained from operator theoretic arguments alone. 

%
%

Like in \cite{BMcR}, we use methods from operator theory and real harmonic analysis to 
obtain \eqref{eq:mainpertest}.
For a self-adjoint operator, say $\Dirb$, the quadratic estimate 
\begin{equation}   
\label{eq:QtQEspectralthm}
  \int_0^\infty \norm{\Qqb_t u}^2\ \frac{dt}{t}\lesssim \norm{u}^2
\end{equation}
is immediate from the spectral theorem coupled with Fubini's theorem. 
Here $\Qqb_t = t\Dirb(\iden + t^2\Dirb^2)^{-1}$ is a 
holomorphic approximation, adapted to the operator $\Dirb$, 
of the projection onto frequencies in a dyadic band around $1/t$.
For the harmonic analyst, the estimate \eqref{eq:QtQEspectralthm}
yields continuity of a wavelet transform, adapted to $\Dirb$, 
and plays the same role in wavelet theory as Plancherel's theorem does in Fourier theory. 
We refer to \cite{Daubechies} by Daubechies in the case $\Qqb_t$ is the projection onto scale $t$ in the multiscale resolution.
These ideas are also central in Littlewood--Paley theory.

Quadratic estimates like \eqref{eq:QtQEspectralthm} are a flexible tool. They can be adapted
to handle non-self-adjoint operators as well as non-commuting operators.
Relevant to this paper is the latter extension, where we want to estimate 
$f(\Dirp)-f(\Dirb)$ as in \eqref{eq:mainpertest}. 
By expressing these operators in terms of resolvents of $\Dirp$ and $\Dirb$ respectively
via the  Dunford functional calculus, such perturbation estimates can be obtained from 
quadratic estimates of the form
\begin{equation}   
\label{eq:QtQEspectralthmpert}
  \int_0^\infty \norm{\Qq_t A \Ppb_t u}^2\ \frac{dt}{t}\lesssim \norm{u}^2.
\end{equation}
Here $\Qq_t$ is like $\Qqb_t$ above but for the operator $\Dirp$, $A$ typically is a bounded multiplication operator, and 
$\Ppb_t = (\iden + t^2\Dirb^2)^{-1}$ should be thought of as a
holomorphic approximation, adapted to the operator $\Dirb$,
to the projections onto frequencies smaller  than $1/t$.

Just like in the non-self-adjoint case in \eqref{eq:QtQEspectralthm}, the estimates
\eqref{eq:QtQEspectralthmpert} are non-trivial and use the specific structure 
of the operators $\Dirp$ and $\Dirb$. When these are differential operators,
allowing non-smooth coefficients, we can use methods from harmonic analysis
to handle  \eqref{eq:QtQEspectralthmpert} essentially as a Carleson embedding theorem.
For operators with simpler  structure than our Dirac operators, it is also possible
to obtain higher order perturbation estimates. In this case the relevant quadratic estimates
look like \eqref{eq:multlinest}.
For our Dirac operators, \eqref{eq:QtQEspectralthmpert} more precisely 
amounts to the two estimates
\begin{align}
&\int_{0}^1 \norm{\Qq_t A_1 \conn (\imath\iden + \Dir)^{-1} \Ppb_t u}^2\ \dtt 
	\lesssim \norm{A_1}_{\infty}^2 \norm{u}^2
\label{eq:A1QE}
\quad \text{and} \ \\
&\int_{0}^1 \norm{t\Pp_t\divv A_2 \Ppb_t u}^2\ \dtt 
	\lesssim \norm{A_2}_{\infty}^2 \norm{u}^2,
\label{eq:A2QE}
\end{align}
which need to be established for $u \in \Lp{2}(\cV)$, where $A_1$ and $A_2$ are  $\Lp{\infty}$ multipliers. 
Through a similarity transformation of $\Dirp$, we can also assume that $\dom(\Dirp)= \dom(\Dirb)$.
Here 
$\Ppb_t = (\iden + t^2\Dirb^2)^{-1}$, $\Pp_t = (\iden + t^2\Dirp^2)^{-1}$, 
$\Qqb_t = t\Dirb(\iden + t^2\Dirb^2)^{-1}$, and $\Qq_t = t\Dirp (\iden + t^2\Dirp^2)^{-1}$.

At a first glance, trying to adapt the proofs in \cite{BMcR} for \eqref{eq:A1QE} and \eqref{eq:A2QE} to the case of manifolds with boundary
seems to be a straightforward exercise. However, closer inspection reveals an interesting dichotomy. 
In \cite{BMcR}, the estimate \eqref{eq:A2QE} was standard and well known to be equivalent to a certain measure being a Carleson measure, 
and the main new work was in establishing \eqref{eq:A1QE}. Here the operator 
$A_1\nabla (\imath\iden+\Dirb)^{-1}$ which is sandwiched between $\Qq_t$ and $\Ppb_t$,
is not a multiplier but also incorporates a singular integral operator $\nabla (\imath\iden+\Dirb)^{-1}$.
To estimate, a  Weitzenb\"ock-type  inequality for $\Dirb$ is needed.
Turning to a manifold with boundary, one sees that \eqref{eq:A1QE} follows as in \cite{BMcR}, \emph{mutatis mutandis}. 
Instead, the presence of boundary forces \eqref{eq:A2QE} to be a non-standard estimate, since new boundary terms 
appear in the absence of boundary conditions for the multiplier $A_2$.
 Indeed, in order for our estimates to be useful, we need to be able to
allow for general $A_2$.  More
precisely, by Stokes' theorem 
$$
  \int_{\cM} \mg(\Pp_t t\divv u,v) d\mu = \int_\Sigma \mg(t \on \rep u, \Pp_t v) d\sigma - \int_{\cM} (u, t\nabla \tilde \Ppb_t v) d\mu.
$$
The second term on the right hand side is bounded by $\norm{u}_{\Lp{2}} \norm{v}_{\Lp{2}}$ by the 
ellipticity and self-adjointness  of $\Dirp$, but clearly the first term  has no such bound.
This means that in \eqref{eq:A2QE}, the operators $\Pp_t t\divv$ are not even bounded, and standard estimates break down. 

An  important  contribution of this paper lies in the new ideas needed to establish \eqref{eq:A2QE}. Here, we observe 
that even though $\Pp_t t\divv$ is unbounded, the operator $\Pp_t t\divv A_2 \Ppb_t$ as a whole is bounded by $\norm{A_2}_{\Lp{\infty}}$
(which is seen from Stokes' theorem and the ellipticity of $\Dirb$).
Building on this observation, we prove \eqref{eq:A2QE} in \S\ref{Sec:Harm2} by adapting, in a non-trivial way, the
standard harmonic analysis proof, usually referred to as a \emph{local $T(1)$} argument.
The inspiration for this analysis comes from \cite{AAH} by Auscher, Axelsson, Hofmann and \cite{AKMC}
by Axelsson, Keith, McIntosh.  To be more precise,  this allows us to 
reduce \eqref{eq:A2QE} for an arbitrary $\Lp{2}$ section
instead  for certain test sections which vanish near the boundary $\Sigma$.
 For this special class of test sections, we are able to adapt
the boundaryless estimates and \eqref{eq:A2QE} becomes standard.  

The remainder of this paper is organised as follows. In \S\ref{Sec:Setup} we state in detail our 
main perturbation estimate in its general form, and in \S\ref{Sec:Dirac}, we show how it is applied to yield the motivating estimate for the 
Atiyah--Singer Dirac operator under perturbation of local boundary conditions.
Then, \S\ref{Sec:Red} contains the proof of Theorem \ref{Thm:Main}, as outlined above.

As aforementioned,  this article is a  sequel to the authors' joint paper \cite{BMcR} with Alan McIntosh.
During our work on this project, McIntosh untimely passed away, leaving us in great sorrow. 
McIntosh's great heritage to mathematics include his widely celebrated unique blend of operator theory and 
harmonic analysis which has lead to breakthroughs like the proof of the Calder\'on conjecture on the $\Lp{2}$ boundedness of the
Cauchy singular integral operator on Lipschitz curves, jointly with Coifman and  Meyer in \cite{CMcM},
and the proof of the Kato square root conjecture on the domain  of the square root of  
elliptic second-order divergence form operators, 
jointly with Auscher, Hofmann, Lacey and Tchamitchian in \cite{AHLMcT}.

The estimates in this paper go back to the multilinear estimates pioneered by McIntosh in connection
with \cite{CMcM}. There, expressions of the form 
\begin{equation}   \label{eq:multlinest}
   \int_0^\infty  \norm{\Qqb_t A_1 \Ppb_t A_2 \Ppb_t A_3 \Ppb_t\cdots A_k \Ppb_t u}_{\Lp{2}}^2 \frac{dt}t
\end{equation}
were bounded by $\norm{u}_{\Lp{2}}^2$.
Formally, the idea is to pass a derivative from $\Qqb_t$, through the general $\Lp{\infty}$ maps $A_i$,
to the rightmost $\Ppb_t$, which becomes $\Qqb_t= t\Dirb \Ppb_t$, and conclude 
the desired estimate by  \eqref{eq:QtQEspectralthm}.
Concretely, this is achieved by harmonic analysis methods and Carleson measures.
The power of this analysis is well known in real-variable harmonic analysis and, in fact,  the necessary 
and much needed algebra  
of $\Ppb_t$ and $\Qqb_t$ operators are in some circles of mathematicians referred to as \emph{McIntoshery}
(or in French \emph{McIntosherie}).

In this paper, we only  employ  the linear case $k=1$ of these multilinear estimates of McIntosh,
leading to first-order perturbation estimates.
Even though  our work is yet another successful example of McIntoshery, we have nevertheless chosen
to not add his name as an author. 
Both authors are former students of  McIntosh, and we know he had as a firm principle for
omitting his name from publications unless he clearly felt that he had contributed to the 
novelties of the paper in a substantial way. Unfortunately, he could not join us this time.
\section*{Acknowledgements}
The first author was supported by the Knut and Alice Wallenberg foundation, KAW 2013.0322 
 postdoctoral program in Mathematics for researchers from outside Sweden as well as SPP2026 from the German Research Foundation (DFG). 
The authors thank Moritz Egert (Paris 11) and Magnus Goffeng (Gothenburg University) 
for useful discussions.  
The authors thank the anonymous referee for a detailed examination of the paper and for useful feedback. 
\section{Setup and statement of main theorem}
\label{Sec:Setup}

\subsection{Manifolds, bundles, and function spaces}
Let $\cM$ be a smooth manifold (possibly noncompact) 
with smooth boundary $\Sigma = \bdy\cM$. 
Throughout, we fix a smooth,  Riemannian metric $\mg$  on $\cM$
and let $\conn$ denote the associated Levi-Civita connection.
We assume that $\mg$ is complete, by which we mean that $(\cM,\mg)$ is complete as a metric space.
By $\intr{\cM}$, we denote the interior $\cM \setminus \bdy\cM$.
The induced volume measure is denoted by $d\mu$ on $\cM$
and $d\sigma$ on $\Sigma$.
Let $\on$ be the unit outward normal vectorfield on $\Sigma$.

The tangent, cotangent bundles are denoted by $\tanb \cM$
and $\cotanb\cM$ respectively, and the rank $(p,q)$-tensor
bundle by $\Tensors[p,q]\cM$. 

For a smooth complex Riemannian bundle $(\cV,\mh)$ on $\cM$, 
let $\Sect(\cV)$ denote the set of measurable sections
and  $\Ck{k,\alpha}(\cV)$ be the set of continuously $k$-differentiable sections
with the $k$-th derivative being $\alpha$-Hölder continuous up to the boundary. 
Note that  when we write $\Ck{k,\alpha}$,
we do not assume $\Ck{k,\alpha}$
with global control of the norm
but rather, only $\Ck{k,\alpha}$
regularity locally. We write $\Ck{k} = \Ck{k,0}$
and $\Ck{\infty}(\cV) = \intersect_{k=1}^\infty \Ck{k}(\cV)$.
Moreover, define
\begin{align*}
&\Ck[c]{k, \alpha}(\cV) = \set{u \in \Ck{k, \alpha}(\cV): \spt u \subset \cM\ \text{compact}}\ \text{and} \\ 
&\Ck[cc]{k, \alpha}(\cV) = \set{u \in \Ck{k, \alpha}(\cV): \spt u \subset \intr{\cM}\ \text{compact}}. 
\end{align*} 
 
Since Lipschitz maps will have special significance, 
we write $\Lips(\cV)$ to denote sections
$\psi \in \Ck{0,1}(\cV)$ with $\norm{\conn{ \psi}}_{\Lp{\infty}(\cV)} < \infty$.  

For $1 \leq p < \infty$,  denote the set of $p$-integrable 
measurable sections with respect to $\mh$ and $\mu$ 
by $\Lp{p}(\cV)$ with  norm $\norm{\xi}_{p}$. The space $\Lp{\infty}(\cV)$
consist of $\xi \in \Sect(\cV)$ such that 
$\modulus{\xi} \leq C$ for some $C > 0$ almost-everywhere on $\cM$.
The norm $\norm{\xi}_{\infty}$ is then the infimum
over $C > 0$ such that this relation holds.
The spaces $\Lp{p}(\cV)$ are Banach spaces and $\Lp{2}(\cV)$
is a Hilbert space with inner product $\inprod{\mdot, \mdot}$. 
The latter space is what we shall be concerned with most 
in this paper and for simplicity of notation, we denote
the norm $\norm{\mdot}_2$ by $\norm{\mdot}$. 
The restricted bundle $\cW = \cV\rest{\Sigma}$ is a smooth, complex Riemannian 
bundle with metric $\mh\rest{\Sigma}$ and $\Lp{p}(\cW)$
spaces are defined similarly on $\Sigma$ with respect to 
the measure $d\sigma$.

Let $\conn$ be a connection on $\cV$ that is compatible with $\mh$.
Then, $\conn$ is a closable operator in   $\Lp{2}(\cV)$  and we
define the Sobolev spaces $\SobH{k}(\cV)$
as the domain of the closure of the operator
$$(\conn, \conn^2, \dots, \conn^k): \Lp{2}\intersect \Ck{\infty}(\cV) \to 
	\Lp{2} \intersect \Ck{\infty}(\oplus_{l=1}^k \Tensors[l,0]\cM \tensor \cV)$$
in $\Lp{2}$. Similarly, we obtain 
boundary Sobolev spaces $\SobH{k}(\cV\rest\Sigma)$  from $\conn \rest{\Sigma}$. 
By compatibility, we have that 
$$ \inprod{\conn u, v} = \inprod{u, -\tr \conn v}$$
 for $u\in \Lp{2} \intersect \Ck{\infty}(\cV)$, $v \in \Lp{2} \intersect \Ck{\infty}(\cotanb \cM \tensor \cV)$ 
and with either   $\spt u \subset \intr{\cM}$ compact or $\spt v \subset \intr{\cM}$ compact.
Thus, we obtain  the divergence operator, defined as $\divv = \adj{\conn_c}$ 
as a densely-defined and closed operator
with domain $\dom(\divv)$ 
from the operator $\conn_c: \Ck[cc]{\infty}(\cV) \to \Ck[cc]{\infty}(\cotanb\cM \tensor \cV)$.

\subsection{Main theorem}
\label{S:MainThm}
In order to phrase the main theorem as in \cite{BMcR}, we require
some assumptions on the manifold.
We say that $(\cM,\mg,\mu)$ has exponential
volume growth if there exists $c_E \geq 1$, $\kappa, c > 0$
such that 
\begin{equation*}
\tag{$\mathrm{\text{E}_{\text{loc}}}$}
\label{Def:Eloc}
0 < \mu(\Ball(x,tr)) \leq ct^\kappa \e^{c_E tr} \mu(\Ball(x,r)) < \infty, 
\end{equation*}
for every $t \geq 1$  and $\mg$-balls $\Ball(x,r)$
of radius $r > 0$ at every $x\in \cM$. 
The manifold $(\cM,\mg)$ satisfies a \emph{local \Poincare inequality}
if there exists $c_P \geq 1$ such that
for all $f \in \SobH{1}(\cM)$, 
\begin{equation*}
\norm{f - f_{B}}_{\Lp{2}(B)} \leq c_P\ \rad(B) \norm{f}_{\SobH{1}(B)}
\tag{$\text{P}_{\text{loc}}$}
\label{Def:Ploc}
\end{equation*}
for all  balls $\Ball$ in $\cM$
such that  the radius $\rad(B) \leq 1$.

We say that $(\cV,\mh)$ satisfies 
\emph{generalised bounded geometry}, or \emph{GBG} for short, 
if there exist $\brad > 0$ and
$C\geq 1 $ such that, for each $x \in \cM$,
there exists a   continuous  local trivialisation 
$\psi_x: \Ball(x,\brad) \times \C^N \to \pi_{\cV}^{-1}(\Ball(x,\brad))$
satisfying
$$
C^{-1} \modulus{\psi_x^{-1}(y)u}_\delta \leq \modulus{u}_{\mh(y)} \leq C \modulus{\psi_x^{-1}(y)u}_\delta,$$
for all $y \in \Ball(x,\brad)$,
where $\delta$ denotes the usual inner product
in $\C^N$ and $\psi_x^{-1}(y) u = \psi_x^{-1}(y,u)$ is the pullback of the
vector $u \in \cV_y$ to $\C^N$ via the local trivialisation $\psi_x$
at $y \in \Ball(x,\brad)$. We call $\brad$ the \emph{GBG
radius}. In typical application, the local trivialisations
will be $\Ck{0,1}$ or smooth.

Letting $\Dir$ and $\Dirp$ be first-order differential operators acting on a bundle $\cV$ over $\cM$ 
and that $\Tr: \SobH{1}(\cV) \to \SobH{\frac{1}{2}}(\cV_\Sigma)$ is the boundary trace map, 
we state the following assumptions adapted to our
setting from \cite{BMcR}, 
\begin{enumerate}[({A}1)]
\item
\label{Hyp:First}
\label{Hyp:BasicF}
\label{Hyp:FinDim}
$\cM$ and $\cV$ are finite dimensional, 
quantified by $\dim \cM < \infty$
and $\dim \cV < \infty$,

\item 
\label{Hyp:ExpVol}
$(\cM,\mg)$ has exponential volume growth quantified by 
$c < \infty$, $c_E < \infty$ and $\kappa < \infty$ in 
\eqref{Def:Eloc},

\item 
\label{Hyp:Poin}
a local \Poincare inequality \eqref{Def:Ploc}
holds on $\cM$ quantified by $c_P < \infty$,
\label{Hyp:BasicL}

\item 
\label{Hyp:TanGBG}
$\cotanb \cM$ has  $\Ck{0,1}$  GBG frames 
$\nu_j$ quantified by $\brad_{\cotanb\cM} > 0$ and
$C_{\cotanb\cM} < \infty$,  with 
$ \modulus{\conn \nu_j} < C_{G,\cotanb\cM}$ 
with $C_{G, \cotanb \cM} < \infty$,

\item 
\label{Hyp:VecGBG}
\label{Hyp:BasGBGL}
$\cV$ has  $\Ck{0,1}$  GBG frames
$e_j$ quantified by $\brad_{\cV} > 0$
and $C_{\cV} < \infty$,   with 
$\modulus{\conn e_j} < C_{G,\cV}$  
with $C_{G,\cV} < \infty$,

\item
\label{Hyp:DirGBG}
$\Dirb$ satisfies
$\modulus{\Dirb e_j} \leq C_{D,\cV}$
with $C_{D,\cV} <\infty$
almost-everywhere inside each GBG
frame $\set{e_j}$,

\item
\label{Hyp:Comm}
We have $\eta \dom(\Dirb) \subset \dom(\Dirb)$
for every  bounded $\eta \in \Ck{\infty}(\cM)$
with $\norm{\conn \eta}_\infty < \infty$, 
and 
$ \comm{ \Dirb, \eta}$ and $\comm{\Dirp, \eta}$ are
pointwise multiplication operators
on almost-every fibre $\cV_x$ with 
 a constant $c_{\Dirb, \Dirp} > 0$ such that 
\begin{equation}
\label{Def:CommConst}
\modulus{ \comm{ \Dirb, \eta} u(x)} \leq c_{\Dirb,\Dirp} \modulus{\conn \eta(x)} \modulus{u(x)}
\end{equation}
for almost-every $x \in \cM$ and the same estimate with $\Dirb$ interchanged
with $\Dirp$,

\item
\label{Hyp:Dom}
$\Dirb$ and $\Dirp$ are self-adjoint operators which are essentially self-adjoint on their restriction to
$$\Ck[c]{\infty}(\cV; \cB) = \set{u \in \Ck[c]{\infty}(\cV): \Tr u \in \cB},$$
where $\cB = \SobH{\frac{1}{2}}(\cE)$
with $\cE \subset \cV\rest{\Sigma}$ a smooth subbundle of $\cV\rest{\Sigma}$, 
and both operators have domain $\dom(\Dirb) = \dom(\Dirp) \subset  \SobH{1}(\cV)$
and with $\const \geq 1$  the smallest  constant satisfying
\begin{equation}
\label{Def:DomConst}
\const^{-1} \norm{u}_{\Dirb} \leq \norm{u}_{\SobH{1}} \leq \const \norm{u}_{\Dirb}
\quad \text{and}\quad
\const^{-1} \norm{u}_{\Dirp} \leq \norm{u}_{\SobH{1}} \leq \const \norm{u}_{\Dirp}
\end{equation}
for all $u \in \dom(\Dirb) = \dom(\Dirp)$
and where $\norm{\mdot}_{\Dirb} = \norm{\Dirb \mdot} + \norm{\mdot}$, 
the operator norm, and

\item 
\label{Hyp:Weitz}

\label{Hyp:Last}
$\Dirb$ satisfies the Riesz-Weitzenb\"ock  condition: 
$\dom(\Dir^2) \subset \SobH{2}(\cV)$ with
\begin{equation}
\label{Def:Weitz}
\norm{\conn^2 u} \leq c_{W} (\norm{\Dir^2 u} + \norm{u})
\end{equation}
for all $u \in \dom(\Dir^2)$ with $c_W < \infty$.
\end{enumerate}

The implicit constants in our perturbation
estimates will be allowed to depend on 
\begin{multline}
\label{Def:MainConst}
\const(\cM,\cV, \Dirb, \Dirp)
	= \max \{\dim \cM, \dim\cV,  c ,  c_E, \kappa,
	c_P,  \brad_{\cotanb\cM} , C_{\cotanb \cM}, C_{G,\cotanb\cM},
\\
 {\brad_{\cV}} , C_{\cV},
C_{G, \cV}, c_{\Dirb}, C_{\Dirb, \cV}, \const,  c_W \} < \infty.
\end{multline}

Our main theorem is the following. 

\begin{theorem}
\label{Thm:Main}
Let $\cM$ be a smooth manifold with smooth 
compact boundary $\Sigma = \bnd \cM$ and let
$\mg$ be a smooth metric  on $\cM$  such that $(\cM, \mg)$ is complete as a metric space. 
Let $(\cV,\mh,\conn)$ be a smooth vector bundle
 over $\cM$  with smooth metric $\mh$ and connection $\conn$
that are compatible.

Let $\Dirb,\ \Dirp$ be two first-order differential 
and assume the hypotheses
\ref{Hyp:First}-\ref{Hyp:Last} on $\cM$, $\cV$, $\Dirb$ 
and $\Dirp$ and that
\begin{equation}
\label{Def:Struc}
\Dirp \psi = \Dirb \psi + A_1 \conn \psi + \divv A_2 \psi  + A_3 \psi,
\end{equation}
holds in a distributional sense for $\psi \in \dom(\Dirb) = \dom(\Dirp)$, where
\begin{equation}
\label{Def:Coeff}
\begin{aligned}
A_1 &\in \Lp{\infty}( \bddlf(\cotanb\cM \tensor \cV, \cV)),  \\
 A_2 &\in \Lp{\infty}\intersect \Lips(\bddlf(\cV, \cotanb\cM \tensor \cV)), \\
A_3 &\in \Lp{\infty}(\bddlf(\cV)),
\end{aligned}
\end{equation}
and let $\norm{A}_\infty = \norm{A_1}_\infty + \norm{A_2}_\infty + \norm{A_3}_\infty$.

Then, for each $\omega \in (0, \pi/2)$
and $\sigma \in (0, \infty]$, 
whenever $f \in \Hol^\infty(\OSec{\omega, \sigma})$,
we have the perturbation estimate
$$\norm{f(\Dirp) - f(\Dirb)}_{\Lp{2}(\cV) \to \Lp{2}(\cV)} \lesssim \norm{f}_{\Lp{\infty}(\Sec{\omega,\sigma})} \norm{A}_\infty,$$
where the implicit constant depends on $\const(\cM,\cV,\Dirb,\Dirp)$. 
\end{theorem}

Here $\OSec{\omega,\sigma}:=\set{x+iy: y^2<\tan^2\omega x^2 + \sigma^2},$
and we say that $f \in \Hol^\infty(\OSec{\omega,\sigma})$
if it is holomorphic on $\OSec{\omega,\sigma}$ and
there exists $C > 0$ such that $\modulus{f(\zeta)} \leq C$.
For a definition of functional calculi $f(\Dir)$ and $f(\Dirp)$
with symbols $f$ bounded and holomorphic, see \S2.3 in \cite{BMcR}.  

\begin{remark}
Self-adjointness of $\Dirb$ and $\Dirp$ in Theorem \ref{Thm:Main} 
\ref{Hyp:Dom} can be relaxed. Indeed, we only use self-adjointness to 
obtain the estimates \eqref{Eq:SA1} 
and \eqref{Eq:SA2}. In the more general situation, i.e. when 
the operator $\Dirb$ or $\Dirp$ is only similar to
a self-adjoint operator with similarity transform $U$,
the constant  $\frac{1}{2} \norm{U}^2 \norm{U^{-1}}^2$ appears in place 
of $\frac{1}{2}$ in \eqref{Eq:SA1} and \eqref{Eq:SA2}, and 
also enters in $\const(\cM, \cV, \Dirb, \Dirp)$. 
\end{remark}

We prove this theorem using real-variable
harmonic analysis methods through the holomorphic
bounded functional calculus in \S\ref{Sec:Red}.
\section{Application to the Atiyah-Singer Dirac operator}
\label{Sec:Dirac}

Throughout this section, in addition to assuming that $(\cM,\mg)$ is a smooth
and complete Riemannian manifold with compact boundary $\Sigma = \bdy \cM$,
we assume that $\cM$ is a \emph{Spin manifold}.

Recall that the exterior algebra
 $\Forms\cM = \oplus_{p=0}^n \Forms[p]\cM$
is a graded algebra, and it is vector-space isomorphic
to the Clifford algebra which we denote  $\Cliff\cM$.
Fix a spin structure $\Prin{\Spin}(\cM)$ and let the associated
Spin bundle be denoted by $\Spinors \cM = \Prin{\Spin} \times_{\eta} \Spinors\R^n$ 
corresponding to the standard complex representation 
 $\eta: \Cliff\R^n \to  \bddlf(\Spinors\R^n)$. 
 Let $\rep: \Sect(\Cliff\cM)  \to \End(\Spinors\cM)$ denote
Clifford multiplication on spinors.  

Let $\SDir$ denote the Atiyah-Singer Dirac operator associated to $\Spinors\cM$, 
given locally in an orthonormal frame $\set{e_k}$ by the 
expression $\SDir \psi = e^k \rep \conn[e_k]\psi$, where $\conn$ is the Spin connection.
 Denoting $\set{\spin{e}_\alpha}$ to be  an induced local orthonormal spin
frame from $\set{e_k}$, the Spin connection takes the local expression   
$\Sconn \spin{e}_\alpha = \conform^2_{E} \rep \spin{e}_\alpha$,
where $\conform^2_{E} =  \frac{1}{2} \sum_{b < a} \conform^a_b \tensor e_b \rep e_a$ 
is the  lifting of the Levi-Civita connection $2$-form to $\Spinors\cM$ 
and   $\conform^a_b$ is the connection $1$-form   in $E=(e_1, \dots, e_n)$.  
The symbol of this operator is $\sym_{\SDir}(\xi)\psi  =  \xi \rep \psi$.
We refer the reader to \cite{LM} by Lawson and Michelsohn,
and \cite{Ginoux} by Ginoux for a more detailed exposition
on spin structures, bundles and their associated operators.  

To define $\SDir$ as a self-adjoint elliptic operator on $\Lp{2}(\Spinors\cM)$
by imposing boundary conditions on $\dom(\SDir)$, we will follow the framework 
developed by Bär and Ballmann in  \cite{BB2} and specialised
to Dirac-type operators in \cite{BB}.   
In particular, by a \emph{local boundary condition} for $\SDir$, we mean a space  
$$\cB = \SobH{\frac{1}{2}}(\cE)\quad \text{with}\quad \cE \subset \Spinors\Sigma = \Spinors \cM \rest{\Sigma},$$ 
where $\cE$ is a smooth subbundle.   
The operator $\SDir$ with boundary condition $\cB$, denoted
$\SDir_{\cB}$, is the operator $\SDir$ with domain 
$$
\dom(\SDir_{\cB}) = \set{ \phi \in \Lp{2}(\Spinors\cM): \SDir \phi \in \Lp{2}(\Spinors\cM)\ \text{and}\ \Tr \phi \in \cB},$$
where $\Tr$ denotes the trace map.
In particular, the choice $\cE = 0$ yield  $\SDir_{\min}$ and $\SDir_{\max} = \close{\SDir_{\SobH{\frac12}{\Spinors\cE}}}$.

Two conditions we require of the local boundary condition $\cB$
are as follows: 

\begin{enumerate}[(i)]
\item Self-adjointness, which by \S3.5 in \cite{BB} occurs
if and only if $\sym_{\SDir}(\on^\flat)$ maps the $\Lp{2}$
closure of $\cB$ onto its orthogonal complement.

\item $\SDir$-ellipticity, which is defined in terms
of a self-adjoint boundary operator $\BDir$ adapted
to $\SDir$ with principal symbol 
$\sym_{\BDir}(\xi) = \sym_{\SDir}(\on^\flat)^{-1} \comp \sym_{\SDir}(\xi)$, 
and for which the operator 
 
$$\pi_{\cB} - \ind{[0, \infty)}(\BDir):\Lp{2}(\Spinors\Sigma) \to \Lp{2}(\Spinors\Sigma)$$ 
is a Fredholm operator. Here, 
$\pi_{\cB}: \Lp{2}(\Spinors\Sigma) \to \cB$ 
is projection induced from the fibrewise orthogonal
projection $\pi_{\cE}: \Spinors\Sigma \to \cE$, and
$\ind{[0,\infty)}(\BDir)$ is the projection
onto the positive spectrum of the operator $\BDir$ (see Theorem 3.15 in \cite{BB}).
This condition yields regularity up to the boundary, 
in the sense that $\SDir u \in \SobH[loc]{k}(\Spinors\cM)$
if and only if $u \in \SobH[loc]{k+1}(\Spinors\cM)$
 whenever $u \in \dom(\SDir_{\cB})$.   
For a compact set $K \subset \cM$, 
the constant $C_K$ such that 
$$ C_K^{-1} \norm{u}_{\SDir_{\cB}^k, K} \leq \norm{u}_{\SobH{k},K} \leq C_K \norm{u}_{\SDir_\cB^k,K}$$
we call  the  \emph{$\SDir$-ellipticity constant of order $k$ in $K$}.   
Here, $\norm{u}_{T, K}^2 = \norm{\ind{K} T u}^2 + \norm{\ind{K} u}^2$. 
See \S7.3-7.4 in \cite{BB2} as well as \S3.5 in \cite{BB}.
\end{enumerate} 

We now state our perturbation result for the Atiyah-Singer
Dirac operator $\SDir_{\cB}$ with local boundary condition
$\cB$. For two local boundary conditions 
$\cB$ and $\tilde{\cB}$,   following \S2 in Chapter IV in \cite{Kato},
we define the \emph{$\Lp{\infty}$-gap}  between the subspaces $\cB$ and $\tilde{\cB}$ as
$$ \gap_\infty(\cB, \tilde{\cB}) = \norm{\gap(\cE_x, \tilde{\cE}_x)}_{\Lp{\infty}(\Sigma)} = \sup_{x\in \Sigma} \modulus{\pi_{\cE}(x) - \pi_{\tilde{\cE}}(x)},$$
where $\pi_{\cE}$ and $\pi_{\tilde{\cE}}$ are the
orthogonal projections from $\Spinors\Sigma$ to $\cE$ and 
$\tilde{\cE}$ respectively.
We let  $\norm{\cB}_{\Lips} = \sup_{x\in \Sigma} \modulus{\conn \pi_{E}(x)}$,  
and similarly for $\tilde{\cB}$. 
For a set $Z \subset \cM$ and $r > 0$, we write $Z_r = \set{x \in \cM: \met_\mg(x, Z) < r}$, 
and $Z_r \disunion Z_r$ to be the double of a neighbourhood $\Sigma$ by pasting along $\Sigma$.

\begin{theorem}
\label{Thm:ASDiracMain}
Let $(\cM,\mg)$ be a smooth, Spin manifold with smooth, compact boundary $\Sigma = \bnd\cM$
 that is complete as a metric space  and suppose that there exists:  
\begin{enumerate}[(i)]
\item  a precompact open neighbourhood $Z$ of $\Sigma$  and $\kappa > 0$ such that $\inj(\cM\setminus Z,\mg) > \kappa$, 
\item $C_R < \infty$ such that $\modulus{ \Ric_\mg} \leq C_R$ 
	and $\modulus{\conn \Ric_\mg}  \leq C_R$ on $\cM \setminus Z$, and
\item  a smooth metric $\mg_Z$ on the double $Z_4 \disunion Z_4$ obtained by pasting along $\Sigma$ and $C_Z < \infty$ and $\kappa_Z > 0$ 
	with $\modulus{\Ric_{\mg_Z}} \leq C_Z$ and $\inj(Z_2 \disunion Z_2, \mg_Z) \geq \kappa_Z$. 
\end{enumerate}  

 Fixing $C_B < \infty$, let  $\cB$ and $\tilde{\cB}$ be two local self-adjoint $\SDir$-elliptic boundary 
which satisfies: 
\begin{enumerate}[(i)]
\item[(iv)] $\norm{\cB}_{\Lips} + \norm{\tilde{\cB}}_{\Lips} \leq C_B$, and
\item[(v)] $\SDir$-ellipticity constants of orders $1$ and $2$  in a given compact neighbourhood $K$ of the boundary.
\end{enumerate}  
Then, for $\omega \in (0, \pi/2)$ and $\sigma > 0$, 
whenever we have $f \in \Hol^\infty(\OSec{\omega, \sigma})$, we have the perturbation
estimate
$$ \norm{f(\SDir_\cB) - f(\SDir_{\tilde{\cB}})}_{\Lp{2}\to \Lp{2}} \lesssim \norm{f}_\infty \gap_\infty(\tilde{\cB}, \cB),$$
where the implicit constant depends on $\dim\cM$ and the constants appearing in (i)-(v).
\end{theorem}

\begin{remark}
The double of a smooth manifold with boundary by pasting along that boundary is again smooth (in terms of the differentiable structure).
However, the canonical reflection of the metric may fail to be smooth across the boundary. 
The existence of a metric $\mg_Z$ satisfying the assumed curvature bounds on $Z_2 \disunion Z_2$ is always guaranteed, but
we have included this in order to quantify the dependence of the constants in the perturbation estimate. See \S\ref{Sec:Dirac1} for more details. 
\end{remark}

\begin{example}[Boundary conditions in even dimensions]
For $\cM$ even dimensional, the Spin bundle
splits $\Spinors \cM = \Spinors^+\cM \oplus^\perp \Spinors^{-} \cM$
 (where $\Spinors^{\pm}\cM$ are the eigenspaces of $u \mapsto \on \rep u$)  
and
$$\SDir = \begin{pmatrix} 0 & \SDir^- \\ \SDir^+ & 0\end{pmatrix},$$
where $\SDir^{\pm}: \Spinors^{\pm}\cM \to \Spinors^{\mp}\cM$.
Again by even dimensionality, $\on: \Spinors^{\pm} \Sigma \to \Spinors^{\mp}\Sigma$.

Let $\B \in \End(\Spinors^+\Sigma)$ smooth and invertible, 
and define 
$$\Spinors_{\B,x} \Sigma = \set{(\psi, \on \rep \B \psi): \psi \in \Spinors^+_x \Sigma}
\ \text{and}\ \Spinors_\B\Sigma = \disunion_{x \in \cM} \Spinors_{\B,x}\Sigma,$$
which is a smooth
sub-bundle of $\Spinors\Sigma$.  The boundary condition   
as considered by Gorokhovsky and  Lesch in \cite{GL}
is then given by $\cB_{\B} = \SobH{\frac{1}{2}}(\Spinors_\B\Sigma)$.

When the boundary condition defining endomorphism $\B$ further satisfies
$\adj{\B(x)} = \B(x)$, then the boundary condition $\cB_{\B}$
is  $\SDir$-elliptic and $\SDir_{B}$ on $\Ck[c]{\infty}(\Spinors\cM; \cB_{\B})$
is essentially self-adjoint.  These facts are a consequence
of Corollary 3.18 in \cite{BB}, which 
guarantees $\SDir$-ellipticity of the boundary condition $\cB_\B$
since $\sym_{\BDir}(\xi)$ interchanges 
$\Spinors_\B \Sigma$ and $\Spinors_\B^\perp \Sigma$
for $0 \neq \xi \in \cotanb_x \Sigma$.
The  essential self-adjointness follows from invoking Theorem 3.11 in \cite{BB},
since $\sym_{\SDir}(\on)$ interchanges $\cB_{\B}$ with its
$\Lp{2}$-orthogonal complement 
$\cB_{\B}^\perp = \set{(\B\on \rep v,  v): v \in \SobH{\frac{1}{2}}(\Spinors^- \Sigma)}$
in $\SobH{\frac{1}{2}}(\Spinors\Sigma)$.
\end{example}

\begin{example} 
As noted in \cite{HMR}, 
Chiral conditions arise from 
an associated Chirality operator  $G \in \Ck{\infty}(\bddlf(\Spinors\cM))$
satisfying:  for all $X \in \Ck{\infty}(\tanb \cM)$ and $\psi,\ \phi \in \Ck{\infty}(\Spinors\cM)$, 
$$ G^2 = \iden,\ \inprod{G \phi, G \psi} = \inprod{\phi, \psi},\ 
\conn[X](G \psi) = G \conn[X]\psi,\ X \rep G \phi = - G X \rep \phi,$$
and the boundary condition is defined via the
projector $\pi_G u = \frac{1}{2}(\iden - \on \rep G)$. 
This is a self-adjoint local elliptic boundary condition which exists in any dimension
(given the map $G$), and has been used in the study of asymptotically 
flat manifolds including black holes. See \S5.2 in \cite{HMR} for more details.
\end{example}

\begin{proof}[Proof of Theorem \ref{Thm:ASDiracMain}]
Without loss of generality, we can assume that $\gap_\infty(\cB, \tilde{\cB)} \leq 1/2$,
as the estimate is trivially true from the spectral theorem for $\gap_\infty(\cB, \tilde{\cB}) > 1/2$.
Note that since the projectors $\pi_{\cE}$ and $\pi_{\tilde{\cE}}$ on $\Spinors\Sigma$ to $\cE$ and $\tilde{\cE}$
respectively are orthogonal, ${\norm{2\pi_{\cE} - \iden}_\infty} = 1$ and so we 
obtain:   
\begin{enumerate}[(i)] 
\item $\norm{\pi_{\cE} - \pi_{\tilde{\cE}}}_\infty \leq \frac{1}{2 \norm{2\pi_{\cE} - \iden}_\infty}$ and
\item $\norm{\conn \pi_{\cE}}_\infty + \norm{\conn \pi_{\tilde{\cE}} }_\infty \leq C_B$
\end{enumerate}
We claim that there exists a
 $\U \in \Lips(\bddlf(\Spinors\cM))$  
with $\norm{\U - \iden}_\infty \leq \gap_\infty(\cB, \tilde{\cB}) \leq \frac{1}{2}$ and $\norm{\conn \U} \lesssim C_B$
such that $\U \cB  = \tilde{\cB}$. 
To see this,  set
$\U_0 = \frac{1}{2}(\iden + (2\pi_{\cE} - \iden)(2\pi_{\tilde{\cE}} - \iden))$
and it is easy to see that $\pi_{\cE} = \U_0^{-1} \pi_{\tilde{\cE}} \U_0$.
Fix $\epsilon > 0$ such that 
$[0, \epsilon) \times \Sigma \cong  N^\epsilon$,
where $N^\epsilon = \set{ x \in \cM: \met(x, \Sigma) < \epsilon}$
and note that $\U_0$ extends to a projection 
$\U'(x) = \U_0(x')$ for $x = (t, x') \in [0, \epsilon) \times \Sigma$. 
Then $\U$ is given by:
\begin{equation*}
\U(x) = \begin{cases} 
	\iden 	&x \not \in N^\epsilon, \\ 
	\cbrac{\iden -  \frac{\met(x,\Sigma)}{\epsilon}} \U'(x) + \frac{\met(x,\Sigma)}{\epsilon} \iden
		&x \in N^\epsilon.
	\end{cases}
\end{equation*}

We verify the hypotheses \ref{Hyp:First}-\ref{Hyp:Last} and invoke
Theorem \ref{Thm:Main} with $\cV = \Spinors\cM$, $\Dir = \SDir_{\cB}$
and $\Dirp = \U^{-1} \SDir_{\tilde{\cB}} \U$ to obtain the estimate
$$ \norm{f(\SDir_{\cB}) - f(\U^{-1} \SDir_{\tilde{\cB}} \U)}_{\Lp{2}\to \Lp{2}} 
	\lesssim \norm{\iden - \U}_\infty \norm{f}_\infty.$$
The passage from this to the required estimate follows from the 
fact that we have $\norm{\iden - \U}_\infty \leq 1/2$ by noting that 
$f(\U^{-1} \SDir_{\tilde{\cB}} \U) = \U^{-1}  f(\SDir_{\tilde{\cB}}) \U$
and that 
$ \norm{f(\SDir_{\tilde{\cB}}) -f(\U^{-1} \SDir_{\tilde{\cB}} \U)}_{\Lp{2} \to \Lp{2}} 
	\lesssim \norm{\iden - \U}_{\infty} \norm{f}_\infty.$

The first hypothesis \ref{Hyp:First} is
immediate and \ref{Hyp:ExpVol} and \ref{Hyp:Poin} are
a consequence of the fact that
the curvature assumptions imply that $\Ric_\mg \geq -C_R$
 (c.f.  Theorem 5.6.4 and 5.6.5 in \cite{SC}).   

The existence of GBG frames satisfying the required bounds
in \ref{Hyp:TanGBG}, \ref{Hyp:VecGBG},
and \ref{Hyp:DirGBG} follow from Proposition \ref{Prop:HarmCoords}, 
 which only depend on $C_R$, $\kappa$, $C_Z$ and $\kappa_Z$.
See \S\ref{Sec:Dirac1}.

Since we assume that $\cB$ is a local boundary condition, 
we have that for every $\eta \in \Lp{\infty} \intersect \Lips(\cM)$, 
the domain inclusion $\eta \dom(\SDir_{\cB}) \subset \dom(\SDir_{\cB})$ holds. 
The commutator estimates follow from the fact
 that 
$$\comm{\SDir,\eta}u = \extd \eta \rep u\quad\text{and}\quad
\comm{\U^{-1} \SDir \U, \eta}u = \U^{-1} \extd \eta \rep  \U u.$$  
This shows \ref{Hyp:Comm}.

The hypothesis \ref{Hyp:Dom} is a consequence of 
Propositions \ref{Prop:StrongEll} and \ref{Prop:DomEq}
since we assume that $\cB$ and  $\tilde{\cB}$ are $\SDir$-elliptic
boundary  conditions.   Note that the constant arising
from these propositions include
the constant $C_{\mathrm{ell},K}$ in the ellipticity estimate
$$C_{\mathrm{ell},K}^{-1} \norm{u}_{\SDir_{\cB},K} \leq \norm{u}_{\SobH{1},K} \leq C_{\mathrm{ell},K} \norm{u}_{\SDir_\cB,K}$$
whenever $u \in \dom(\SDir_{\cB})$.
The corresponding constant in the region $\cM \setminus K$
depends on the geometric bounds (i)-(iii).   
 In addition to these constants for $\SDir_{\cB}$, 
the corresponding estimate for the operator $\SDir_{\tilde{\cB}}$
includes the constant $C_B$. See \S\ref{Sec:Dirac2} for details. 

The remaining hypothesis is the Riesz-Weitzenböck hypothesis \ref{Hyp:Weitz}.
This is proved similar to Proposition \ref{Prop:StrongEll},
using the compact set $K$ and $K_{\frac{1}{2}}$ near
the boundary, along with the smooth cutoff $f$ 
as they appear in the proof of this proposition.
The estimate $ \norm{\conn^2 (fu)} \lesssim \norm{\SDir^2_{\cB} u} + \norm{u}$
is obtained by arguing as in Proposition 3.18 in \cite{BMcR}
via the cover provided by Lemma \ref{Lem:GradBddCover}, 
and the remaining estimate
$\norm{\conn^2((1 - f)u)} \leq \tilde{C}_{\mathrm{ell},K} (\norm{\SDir^2_{\cB} u} + \norm{u})$
is due to the boundary regularity result, Theorem 7.17 in \cite{BB2}.
Here,  ellipticity constant $\tilde{C}_{\mathrm{ell},K}$ is
the constant 
$$\tilde{C}_{\mathrm{ell},K}^{-1} \norm{u}_{\SDir_{\cB}^k,K} \leq \norm{u}_{\SobH{k},K} \leq \tilde{C}_{\mathrm{ell},K} \norm{u}_{\SDir_\cB^k,K}$$
whenever $u \in \dom(\SDir_{\cB}^k)$ for $k = 1,2$.
The constant for the estimate in the region $\cM \setminus K$ depend on 
the constants in (i)-(iii). 

Lastly, the decomposition of the  operator 
$\SDir_{\tilde{\cB}} - \SDir_{\cB} = A_1 \conn + \divv A_2   + A_3$
 distributionally  proved in Proposition \ref{Prop:OpDecomp}. 
See \S\ref{Sec:Dirac3} for details.
\end{proof} 

Throughout the remainder of this section, we assume the hypothesis of Theorem \ref{Thm:ASDiracMain}.

\subsection{Geometric bounds in the presence of boundary}
\label{Sec:Dirac1}
\newcommand{\cMt}{\widetilde{\cM}}
\newcommand{\mge}{\mg_{\rm{ext}}} 

The way in which we prove Theorem \ref{Thm:ASDiracMain} is via Theorem \ref{Thm:Main}, which
requires us to prove that under the geometric assumptions we make, the bundle $\Spinors\cM$
satisfies generalised bounded geometry and the  first and second metric derivatives in 
each trivialisation are bounded.

We do this by considering  the double of the manifold $\cMt = \cM \disunion \cM$,
which is obtained  
by taking two copies of $\cM$ and pasting along the boundary $\Sigma$ to obtain a manifold without boundary.
Since the boundary  is smooth, this manifold is again smooth 
(in a differential topology sense, see Theorem 9.29 in \cite{Lee}).
By reflection, we obtain an extension $\mg_{\rm{ext}}$ of the metric $\mg$ to the whole of $\cMt$.
This metric is guaranteed to be continuous everywhere and smooth on $\cMt \setminus \Sigma$, 
but in general, without imposing additional restrictions on the boundary, it will not be smooth. 
However, as we illustrate in the following lemma, we are able to
construct a smooth metric sufficiently close to $\mge$ that suffices to obtain 
the bounds we desire for $(\cM,\mg)$.

\begin{lemma}
\label{Lem:NearMet}
There exists a smooth complete metric  $\mgt$ 
on $\cMt$ with $ G \geq 1$ dependent on $\mg_Z$ and $\mg$ satisfying
	$$ G^{-1} \modulus{u}_{\mgt}  \leq \modulus{u}_{\mg_{\rm{ext}}} \leq G \modulus{u}_{\mgt}$$
and for which there exists:
\begin{enumerate}[(i)]
\item $\tilde{\kappa} > 0$ such that $\inj(\cMt,\mgt) >\tilde{\kappa}$,

\item $\tilde{C}_{R} < \infty $ such that $\modulus{\Ric_{\mgt}} \leq \tilde{C}_R$ and
	$\modulus{\conn \Ric_{\mgt}} < \tilde{C}_R$,

\item  a compact set $\cP$ with $\interior{\cP} \neq\emptyset$ and $\Sigma \subset \cP$
	such that $\mg_{\rm{ext}}  = \mgt$ on $\cMt \setminus \cP$.  
\end{enumerate}

The constants $\tilde{\kappa}$, $\tilde{C}_R$ and depend on the original geometric bounds $\kappa$, $C_R$, $\kappa_Z$, $C_Z$. 
\end{lemma}
\begin{proof}

Take $Z$ from the hypothesis of Theorem \ref{Thm:ASDiracMain} and let $\cP = \close{Z \disunion Z}$.
By hypothesis, since $Z$ is precompact, we get that $\cP$ is compact.
As a consequence, if $\set{x_n}$ is a Cauchy sequence in $\cP$, then it converges to some point
and if $\set{x_n}$ is Cauchy in  $\cMt \setminus \interior{\cP}$, then it converges to some point in $\cMt \setminus \interior{\cP}$ by the metric completeness of $\mg$. 
This establishes that $\mge$ is metric complete.

Next, let $\psi \in \Ck{\infty}(\cMt)$ be such that $\psi = 1$ on $\cMt \setminus \interior{\cP}$
and $\psi = 0$ on $\cP_{\frac{3}{2}} = \set{x \in \cMt: \met_{\mge}(x,\cP) \leq \frac{3}{2}}$.
Since $\cP_\epsilon$ is compact by construction, by the smoothness of the differentiable
structure of $\cMt$, there exists $G \geq 1$ such that $\mge$ and $\mg_Z$ are $G$-close on $\cP_2$.
Define $\mgt = \psi \mge + (1 - \psi) \mg_Z$
and since $\mge = \mgt$ away from $\cP$, this
shows that the quasi-isometry with constant $G$ between $\mge$ and $\mgt$ and also establishes (iii).

Since $\mg_Z$ satisfies a lower bound on injectivity radius on $Z_2 \disunion Z_2$ as well as a Ricci curvature
bound on this set, and since $\mg$ satisfies similar bounds on $Z$, by construction of the metric $\mgt$,
we obtain (i) and (ii)  with the dependency as stated in the conclusion.
\end{proof}

Now, using this we can prove the main proposition that we require
to prove the geometric bounds needed to prove Theorem \ref{Thm:ASDiracMain}.
\begin{proposition}
\label{Prop:HarmCoords}
There exist $r_H > 0$ and a constant $1 \leq C < \infty$ depending on $\kappa$, $C_R$, $\kappa_Z$ and $C_Z$ such that
at each $x \in \cM$, $\psi_x: B(x,r_H) \to \R^n$ corresponds to a coordinate system
and inside that coordinate system with coordinate basis $\set{\partial_j}$
satisfying:  
$$C^{-1}\modulus{u}_{\pullb{\psi}_x\delta(y)} \leq \modulus{u}_{\mg(y)} \leq C \modulus{u}_{\pullb{\psi}_x\delta(y)},\quad
\modulus{\partial_k \mg_{ij}(y)} \leq C,\quad\text{and}\quad \modulus{\partial_k \partial_l \mg_{ij}(y)} \leq C,$$
for all $y \in B(x,r_H)$ and where $\delta$ is the Euclidean metric.
\end{proposition}
\begin{proof}
Utilising the metric $\mgt$ given by Lemma \ref{Lem:NearMet}, we
apply Theorem 1.2 in \cite{Hebey} to obtain $\Ck{2,\alpha}$-harmonic coordinates 
for the manifold $(\cMt, \mgt)$ with radius $\widetilde{r_H}$. We obtain the same conclusions
for $(\cM,\mgt\rest{\cM})$ as it is obtained via the subspace topology on 
$\cMt$. The balls $B_{\mg}$ and $B_{\mgt}$ are contained
within the factor $G$ given in the lemma, and away from the compact region 
$\cP$ defined in the lemma, we have that $B_{\mg} = B_{\mgt}$. So, it
suffices to set $r_H = \widetilde{r_H}/G$.
On the region $\cMt \setminus \cP$, we have $\Ck{2,\alpha}$ control of the metric $\mgt$
and outside of this region, by compactness, we obtain control of as many derivatives
of the metric as we like. By taking maximums of the constants appearing
in the regions $\cMt \setminus \cP$ and $\cP$ , we obtain the constant $C$ in the conclusion of this
proposition. 
\end{proof}

\subsection{The domains of the operators}
\label{Sec:Dirac2}
To invoke Theorem \ref{Thm:Main}, we need to establish 
Sobolev $\SobH{1}$ regularity for the operators $\SDir_{\cB}$ and $\SDir_{\tilde{\cB}}$. To this end,
we begin with the following covering lemma. 

\begin{lemma} 
\label{Lem:GradBddCover}
There exists $C_H < \infty$, $M > 0$ and a sequence of points $x_i$ and
a smooth partition of unity $\set{\eta_i}$ for $\cM$ that is uniformly locally finite
and subordinate to $\set{B(x_i,r_H)}$ satisfying:
\begin{enumerate}[(i)]
\item $ \sum_{i} \modulus{\conn^j \eta_i} \leq C_H$ for $j = 0, \dots, 3$, and 
\item $1 \leq M \sum_{i} \eta_i^2$.
\end{enumerate}
The $r_H > 0$ here is the harmonic radius guaranteed in Proposition \ref{Prop:HarmCoords}.
\end{lemma} 
\begin{proof}
Take
the double of the manifold and the smooth metric given by Lemma \ref{Lem:NearMet}. 
Then, by Lemma 1.1 in \cite{Hebey}, on fixing $\rho > 0$ we find a sequence of points
$x_i \in \cMt$ such that (i) $\set{\tilde{B}(x_i, r)}$ is a uniformly locally finite cover of $\cMt$ 
for all $r \geq \rho$ and (ii) $\tilde{B}(x_i, \rho /2) \intersect \tilde{B}(x_j,\rho/2) = \emptyset$
for all $i \neq j$.
This relies purely on a measure counting argument since $\mgt$ induces
a measure satisfying exponential volume growth \eqref{Def:Eloc} 
by the Ricci curvature lower bounds. Since
 $\mgt$ is $G$-close to $\mge$, the same is true for the metric $\mge$, 
which is the metric guaranteed to be continuous obtained by reflection of $\mg$ on $\cM$ 
across $\Sigma$ to the double $\cMt$. Thus, a cover
satisfying (i) and (ii) exists on $\cMt$ replacing $\mgt$ balls $\tilde{B}$ 
with $\mge$ balls $B^{\rm{ext}}$.

Now, let $r_H$ denote the radius obtained from Proposition \ref{Prop:HarmCoords},
and set $\rho = r_H/16$. 
Let $\set{x_i^\cM} \subset \interior{\cM}$ such that $\met_{\mg}(x_i^\cM, \Sigma) > r_H/16$.
Then $\set{x_i^\cM} \subset \cM \setminus Z'$, where 
$Z' = \set{x \in \cM: \met_{\mg}(x, \Sigma) \leq r_H/16}$.
Since $\Sigma$ is compact, so is $Z'$ and hence, there exists
a finite number of points $\set{x_j^{Z'}}_{j=1}^K$ such that 
$Z' \subset \union_{j = 1}^K B(x_j^{Z'}, r_H/16)$.
Then, the collection of points
$ \set{\bar{x}_i} = \set{x_j^\cM, x_k^{Z'}}$ satisfies:
$ \cM = \union_{i} B(\bar{x}_i, r_H/16)$
with $\set{B(\bar{x}_i, r_H/16)}$ uniformly locally finite.

Inside each $B(\bar{x}_i, r_H/16)$ we have $\Ck{2,\alpha}$ control 
of the metric, and therefore, the partition of unity $\set{\eta_j}$  
with the gradient bound in the conclusion is
obtained by proceeding as in the the proof of Proposition 3.2 in \cite{Hebey}.
\end{proof}

With this lemma, we prove the following. 
\begin{proposition}
\label{Prop:StrongEll}
 
The embedding $\dom(\SDir_{\cB}) \embed \SobH{1}(\cV)$ holds
along with the ellipticity estimate $\norm{u}_{\SDir_{\cB}} \simeq \norm{u}_{\SobH{1}}$
for all $u \in \dom(\SDir_{\cB})$.    
\end{proposition}
\begin{proof}
Let $K$ be a compact neighbourhood of $\Sigma$ assumed in (v) of Theorem \ref{Thm:ASDiracMain}
and let $f: \cM \to [0,1]$ be smooth
with $f = 1$ on $\cM \setminus \interior{K}$ and
$f = 0$ on an open subset $\tilde{K} \subset \interior{K}$ with $\Sigma \subset \tilde{K}$.
Let $u \in \dom(\SDir_{\cB})$ and we show that
$\norm{\conn(f u)} + \norm{\conn((1 - f)u)}  \lesssim \norm{\SDir_{\cB} u } + \norm{u}.$
Using the cover guaranteed by Lemma \ref{Lem:GradBddCover}, we obtain that 
$$\norm{\conn(f u)} \lesssim \norm{\SDir_{\cB}(fu)} + \norm{fu} \lesssim \norm{\SDir_{\cB} u} + \norm{u},$$
where the first inequality is from running 
the exact same argument as Proposition 3.6 in \cite{BMcR}
and the second inequality is from the fact that
$\spt \conn f \subset K$ and hence bounded.
For the remaining inequality, we note that
since the boundary condition $\cB$ 
is $\SDir$-elliptic, Theorem 7.17 in \cite{BB2}
gives us that $u \in \SobH[loc]{k+1}(\Spinors\cM) \iff \SDir_{\cB} u \in \SobH[loc]{k}(\Spinors\cM)$
whenever $u \in \dom(\SDir_{\cB})$.
Choosing $k = 0$, and the fact that $\spt (1 -f)u \subset K$, we get that
$$ \norm{\conn ((1 - f) u)} \leq C_{\mathrm{ell},K} (\norm{\SDir_{\cB}((1 - f) u)} +  \norm{(1 - f) u})
\lesssim \norm{\SDir_{\cB} u } + \norm{u}.$$ 
where $C_{\mathrm{ell},K} < \infty$ is a constant that depends on $K$.

The estimate $\norm{u}_{\SDir_{\cB}} \lesssim \norm{u}_{\SobH{1}(\cV)}$ 
for $u \in \dom(\SDir_{\cB})$ follows from
the pointwise estimate $\modulus{\SDir u} \lesssim \modulus{\conn u}$
(c.f. Proposition 3.6 in \cite{BMcR}).  
\end{proof}

Using this proposition, we prove the following.
\begin{proposition}
\label{Prop:DomEq}
The equality $\dom(\SDir_{\cB}) = \dom(\SDir_{\tilde{\cB}} \U)$ holds.
\end{proposition}
\begin{proof}
On fixing $\phi \in \Ck[c]{\infty}(\Spinors\cM)$, we compute 
at a point $x \in \cM$ with a frame satisfying $\conn[e_i] e_j (x) = 0$: 
$$\SDir(\U \phi) 
	= e^i \rep \conn[e_i] (\phi^\alpha \U \spin{e}_\alpha)
	= (e_i \phi^\alpha) e^i \rep \U \spin{e}_\alpha + \phi^\alpha (e_i \U^\beta_\alpha) \spin{e}_\beta,$$
from which it follows directly that
$ \modulus{\SDir(\U \phi)}^2 \leq \modulus{\U}^2 \modulus{\conn \phi}^2 + \modulus{\conn \U}^2 \modulus{\phi}^2.$
Now, for $\phi \in \dom(\SDir_{\tilde{\cB}})$, we have from Theorem 3.10 in \cite{BB} that
there is a sequence $\phi_n \in \Ck[c]{\infty}(\Spinors\cM; \tilde{\cB})$ such that 
$\phi_n \to \phi$ in the graph norm of $\SDir_{\tilde{\cB}}$. Moreover, $\U \phi_n \in \Ck[c]{0,1}(\Spinors\cM, {\cB}) \subset \dom(\SDir_{\cB})$
and by Proposition \ref{Prop:StrongEll}, 
$\norm{\conn(\phi_n - \phi)} \to 0$. Hence,
combining this with our pointwise estimate and integrating,
we obtain that
$$ \norm{\SDir(\U \phi_n - \U \phi_m)} \lesssim \norm{\U}_\infty \norm{\conn(\phi_n - \phi_m)} + \norm{\conn \U}_\infty \norm{\phi_n - \phi_m}
\to 0$$
as $m, n \to \infty$. By the closedness of $\SDir_{\cB}$, we have that $\U \phi \in \dom(\SDir_{\cB})$. 
The reverse containment is obtained similarly. 
\end{proof}

\subsection{Decomposition of the  difference of operators}
\label{Sec:Dirac3}

A crucial assumption in Theorem \ref{Thm:Main} is to be able
to write the difference of our operators $\SDir_{\cB}$ and $\U^{-1}\SDir_{\tilde{\cB}} \U$
as
$$\SDir_{\cB} - \U^{-1} \SDir_{\tilde{\cB}} \U = A_1 \conn + \divv A_2 + A_3$$ 
with $\norm{A_i}_\infty$ controlled by $\norm{\U - \iden}_\infty$.

Our computations here are similar to those in \S3 of \cite{BMcR}, 
with the key observation being that the last term in 
Lemma \ref{Lem:Decomp1} cannot be used as $A_3$, since it
would yield only a bound $\norm{A_3}_\infty \lesssim 1$
and not $\norm{A_3}_\infty \lesssim \norm{\U - \iden}_\infty$.
Instead, we proceed via an application of the product rule for derivatives as in 
Lemma \ref{Lem:Decomp2}. 
 
Throughout this subsection, 
unless otherwise stated, we fix an open set $\Omega \subset \interior{\cM}$
and let $\set{e_i}$ and  $\set{\spin{e}_\alpha}$ be
orthonormal frames for $\tanb\cM$ and $\Spinors\cM$ respectively inside $\Omega$.

\begin{lemma}
\label{Lem:Decomp1}
For $\phi \in \Ck{\infty}(\Spinors \cM)$ we have the following
pointwise equality almost-everywhere inside $\Omega$: 
$$ 
(\SDir - \U^{-1} \SDir \U)\phi =  X \conn \phi +  Z^\Omega \phi + \phi^\alpha  (\conn[e_j]\U^\beta_\alpha) \U^{-1} e^j \rep \spin{e}_\beta$$
with $X: \Sect(\cotanb\cM \tensor \Spinors\cM) \to \Sect(\Spinors\cM)$
and $Z^\Omega: \Sect(\Spinors\Omega) \to \Sect(\Spinors\Omega)$ with almost-everywhere pointwise estimates 
$$\modulus{X} \lesssim \norm{\iden - \U}_{\infty}\quad{\text{and}}\quad
\modulus{Z^\Omega} \lesssim \norm{\iden - \U}_{\infty},$$
where the implicit constants depends on the constants in Theorem \ref{Thm:ASDiracMain}. 
\end{lemma}
\begin{proof}
A direction calculation yields that
$$
\SDir(\U \phi) = (\conn[e_j] \phi^\alpha) e^j \rep \U \spin{e}_\alpha 
	+ \phi^\alpha (\conn[e_j] \U^\beta_\alpha) e^j \rep \spin{e}_\beta
	+ \phi^\alpha \U^\beta_\alpha e^j \rep \conn[e_j]\spin{e}_\beta.$$
Since the term $\conn[e_j] \spin{e}_\beta = \conform^2_E(e_j) \rep\spin{e}_\beta$, 
multiplying this expression by $\U^{-1}$ on the left, and then subtracting
it from the expression for $\SDir \phi$,
we obtain that 
\begin{multline*}
(\SDir - \U^{-1} \SDir \U)\phi 
	= \conn[e_j]\phi^\alpha (e^j \rep \spin{e}_\alpha - \U^{-1} e^j \rep \U \spin{e}_\alpha) \\
	+ (e^j \rep \conform^2_E(e_j) - \U^{-1} e^j \rep  \conform^2_E(e_j)\U) \rep \phi 
	+ \phi^\alpha (\conn[e_j] \U^\beta_\alpha) \U^{-1} e^j \rep \spin{e}_\beta.
\end{multline*}

To obtain a bound on the first expression to the right of this,
we note that
$$
e^j \rep \spin{e}_\alpha - \U^{-1} e^j \rep \U \spin{e}_\alpha
	= \e^j \rep (\iden - \U) \spin{e}_\alpha + (\iden - \U^{-1}) \e^j \rep \U \spin{e}_\alpha,$$
and we can write 
$$
\conn[e_j] \phi^\alpha \e^j \rep (\iden - \U) \spin{e}_\alpha 
	= X_1\conn \phi - \phi^\alpha e^j  \rep (\iden - \U) \conform^2_E(e_j) \rep \spin{e}_\alpha,$$
where $X_1( \psi^\alpha_k e^k \tensor \spin{e}_\alpha) = \psi^\alpha_k e^k \rep (\iden - \U) \spin{e}_\alpha.$
Now, similarly, writing 
$X_2(\psi^\alpha_k e^k \tensor \spin{e}_\alpha) = \psi^\alpha_k (\iden - \U^{-1}) e^j \rep \U \spin{e}_\alpha$,
we obtain that
$$
\conn[e_j] \phi^\alpha (\iden - \U^{-1}) \e^j \rep \U \spin{e}_\alpha 
	= X_2 \conn \phi - \phi^\alpha (\iden - \U^{-1}) e^j  \rep \U \conform^2_E(e_j) \rep \spin{e}_\alpha.$$
Letting $X = X_1 + X_2$, we obtain that
\begin{multline*}
\conn[e_j]\phi^\alpha (e^j \rep \spin{e}_\alpha - \U^{-1} e^j \rep \U \spin{e}_\alpha)
	= X \conn \phi  \\
		- e^j \rep (\iden - \U) \conform^2_E(e_j) \rep \phi
		- (\iden - \U^{-1}) e^j \rep \U \conform^2_E(e_j) \rep \phi.
\end{multline*}

Now, note that
$$e^j \rep \conform^2_E(e_j) - \U^{-1} e^j \rep \conform^2_E(e_j)\U 
	= e^j \rep \conform^2_E(e_j)( \iden - \U) + (\iden - \U^{-1}) e^j \rep \conform^2_E(e_j) \U,$$
and on setting
\begin{multline*}Z^\Omega = 
\e^j \rep \conform^2_E(e_j)( \iden - \U) + (\iden - \U^{-1}) e^j \rep  \conform^2_E(e_j) \U \\ 
	- e^j \rep (\iden - \U) \conform^2_E(e_j) - (\iden - \U^{-1}) e^j \rep \U \conform^2_E(e_j),
\end{multline*}
we obtain the conclusion. 
\end{proof}

This lemma illustrates that the main term to analyse 
is the last term given by $\phi^\alpha  (\conn[e_j]\U^\beta_\alpha) \U^{-1} e^j \rep \spin{e}_\beta$.
This is the content of the following lemma.
\begin{lemma}
\label{Lem:Decomp2}
For $\phi \in \Ck{\infty}(\Spinors\cM)$, we have the following decomposition
pointwise almost-everywhere inside $\Omega$:
$$\phi^\alpha  (\conn[e_j]\U^\beta_\alpha) \U^{-1} e^j \rep \spin{e}_\beta
	=  L^\Omega \conn \phi + \divv M^\Omega \phi + N^\Omega \phi.$$
The coefficients satisfy the estimates
$$ \norm{L^\Omega}_\infty + \norm{M^\Omega}_\infty + \norm{N^\Omega}_\infty \lesssim \norm{\iden - \U}_\infty
\ \text{and}\ \norm{\conn M^\Omega}_\infty \lesssim 1,$$
where the implicit constants depend on the constants listed in Theorem \ref{Thm:ASDiracMain}.
\end{lemma}
\begin{proof}
First note that on letting $\epsilon^\beta_\alpha = \delta^\beta_\alpha - \U^\beta_\alpha$, 
we have
$
\phi^\alpha  (\conn[e_j]\U^\beta_\alpha) \U^{-1} e^j \rep \spin{e}_\beta
	= -\phi^\alpha (\conn[e_j]\epsilon^\beta_\alpha) \U^{-1} e^j \rep \spin{e}_\beta.$
Let $M^\Omega: \Sect(\Spinors\Omega) \to \Sect(\cotanb\Omega \tensor \Spinors\Omega)$
written inside $\Omega$ as
$$ M^\Omega \psi = \phi^\alpha M_{\alpha, k}^\theta\ e^k \tensor \U^{-1} e^k \rep \spin{e}_\theta$$
with the coefficients to be determined later.
Note that:
\begin{multline*}
\conn(M^\Omega \phi) = M^\theta_{\alpha,k} (\conn[e_j]\phi^\alpha) e^j \tensor e^k \tensor \U^{-1} e^k \rep \spin{e}_\theta \\
	+ \phi^\alpha \conn[e_j] (M^\theta_{\alpha,k}) e^j \tensor e^k \tensor \U^{-1} e^k \rep \spin{e}_\theta 
	+ \phi^\alpha M^\theta_{\alpha, k} e^j \tensor \conn[e_j](e^k \tensor \U^{-1} e^k \rep \spin{e}_\theta).
\end{multline*}
On taking the trace, and rearranging the equation,
\begin{multline*} 
\phi^\alpha \conn[e_j](M^\theta_{\alpha,j}) \U^{-1} e^j \rep \spin{e}_\theta  = \tr_\mg \conn(M^\Omega \phi) \\
	- (\conn[e_j] \phi^\alpha) M^\theta_{\alpha,j} \U^{-1} e^j \rep \spin{e}_\theta   
	- \phi^\alpha M^\theta_{\alpha, k} \tr( e^j \tensor \conn[e_j]( e^k \tensor \U^{-1} e^k \rep \spin{e}_\theta)).
\end{multline*}
So set $M^\theta_{\alpha, k} = \epsilon^\theta_\alpha$, which gives us
an expression for $\phi^\alpha (\conn[e_j]\epsilon^\beta_\alpha) \U^{-1} e^j \rep \spin{e}_\beta$.

It remains to show that the remaining terms in this expression can be decomposed
to $L^\Omega \conn \phi + N \phi$. Let $L^\Omega (e^j \tensor \spin{e}_\alpha)  = \U^{-1} e^j \rep M^\theta_{\alpha,j} \spin{e}_\theta$, 
then we have that 
$$ \tr_{\mg}(\phi^\alpha \conn[e_j] (M^\theta_{\alpha,k}) e^j \tensor e^k \tensor \U^{-1} e^k \rep \spin{e}_\theta)
	= L^\Omega \conn \phi - \U^{-1} e^j \rep M^\Omega \omega^2_E(e_j) \rep \phi.$$
Absorbing the error term in this computation along with the remaining term from 
the former expression, we can set 
$$N^\Omega \phi = - \phi^\alpha \epsilon^\theta_\alpha \tr( e^j \tensor \conn[e_j]( e^k \tensor \U^{-1} e^k \rep \spin{e}_\theta))
	 - \U^{-1} e^j \rep M^\Omega \omega^2_E(e_j) \rep \phi.$$
The estimates in the conclusion for $L^\Omega, M^\Omega, N^\Omega$
and $\conn M^\Omega$ follows from the definitions of these maps.
\end{proof} 

Using these two lemmata, arguing in a similar way to 
Proposition 3.16 in \cite{BMcR}, we obtain the following
decomposition globally on $\interior{\cM}$. 
\begin{proposition}
\label{Prop:OpDecomp}
We have that:
$$ (\SDir_{\cB}  - \U^{-1} \SDir_{\tilde{\cB}} \U) \phi = A_1 \conn \phi + \divv A_2 \phi + A_3 \phi$$
distributionally for all $\phi \in \dom(\SDir_{\cB})$ where the coefficients $A_i$ satisfy:
\begin{align*}
A_1 &\in \Lp{\infty}(\bddlf(\cotanb\cM \tensor \Spinors\cM)), \\ 
A_2 &\in \Lp{\infty} \intersect \Lips(\bddlf(\Spinors\cM, \cotanb \cM \tensor \Spinors\cM)), \\
A_3 &\in \Lp{\infty}(\bddlf(\Spinors\cM)),
\end{align*}
with $\norm{A_1}_\infty + \norm{A_2}_\infty + \norm{A_3}_\infty \lesssim \norm{\iden - \U}_\infty$.
The implicit constants depend on the  constants listed  in Theorem \ref{Thm:ASDiracMain}. 
\end{proposition}
\begin{proof}
Following the proof of Proposition 3.16 in \cite{BMcR}, it suffices
to show that there exists a cover $\set{B_j}$ of balls with 
a fixed radius $r > 0$ with orthonormal frames $\set{e_{j,l}}$ inside
$B_j$, and a Lipschitz partition of unity $\set{\eta_j}$
subordinate to $\set{B_j}$ satisfying: $\modulus{\conn e_{j,l}} \leq C_1$
and $\modulus{\conn \eta_j} \leq C_2$, where $C_1$ and $C_2$
are finite constants independent of $j$ and $l$.
The covering with the gradient bound on the partition of 
unity is given in Lemma \ref{Lem:GradBddCover}
and  the uniform control of $\modulus{\conn e_{i,k} }\leq C_1$
is a consequence of the fact that 
each $B_j$ corresponds to a ball in which 
we have $\Ck{2,\alpha}$ uniform control of the metric.
Then, as in Proposition 3.16 in \cite{BMcR}, using Lemma \ref{Lem:Decomp1} and Lemma \ref{Lem:Decomp2}, we set 
\begin{align*}
A_1 \phi &= X\phi + \sum_{j}  L^{B_j}\eta_j \phi \\
A_2 \phi &= \sum_j M^{B_j}  \eta_j \phi \\
A_3 \phi &= \sum_j (N^{B_j} + Z^{B_j}) \eta_j \phi - \sum_j \tr(  \conn \eta_j \tensor \phi).
\end{align*}
It is readily verified that this yields the desired decomposition.
\end{proof}

%
%
\section{Operator theory and harmonic analysis}
\label{Sec:Red}

Throughout this section, we assume
the hypothesis of Theorem \ref{Thm:Main}. 
Moreover, we assume that the reader is familiar
with the holomorphic functional calculus 
via the Riesz-Dunford integral and 
how to estimate functional calculus of non-smooth
operators with harmonic analysis.
A brief description of this framework is included
in \S2.1 in \cite{BMcR}, but \cite{ADMc} is a
more detailed reference.

For $t > 0$, define the operators 
\begin{align*}
&\Rrb_t = \frac{1}{\iden + \imath t \Dirb},\ 
\Rr_t = \frac{1}{\iden + \imath t \Dirp},\\
&\Ppb_t = \frac{1}{\iden + t^2\Dirb^2},\  
\Pp_t = \frac{1}{\iden + t^2\Dirp^2},\\ 
&\Qqb_t = t\Dirb \Ppb_t,
\ \text{and}\ 
\Qq_t = t\Dirp \Pp_t.
\end{align*}
 Due to self-adjointness, we have the bounds  
\begin{equation}
\label{Eq:SA1}
\int_{0}^\infty \norm{\Qq_t u}^2\ \dtt \leq \frac{1}{2} \norm{u}^2
\quad\text{and}\quad
\int_{0}^\infty \norm{\Qqb_t u}^2\ \dtt \leq \frac{1}{2} \norm{u}^2,
\end{equation}
and 
\begin{equation}
\label{Eq:SA2}
\sup_{t} \norm{\Rrb_t},\ 
\sup_{t} \norm{\Rr_t},\ 
\sup_{t} \norm{\Ppb_t},\ 
\sup_{t} \norm{\Pp_t},\ 
\sup_{t} \norm{\Qqb_t},\ 
\sup_{t} \norm{\Qq_t} \leq \frac{1}{2}.
\end{equation}
Each of these operators are also self-adjoint.

We note the identities  
\begin{equation}
\label{Eqn:Rrt}
\Rr_t = \Pp_t - \imath \Qq_t
\quad \text{and}\quad
 \Rrb_t = \Ppb_t - \imath \Qqb_t,
\end{equation}
  as well as    
\begin{equation}
\label{Prop:Paraprod}
\Rr_t - \Rrb_t = \Rr_t[\imath t(\Dirb - \Dirp)]\Rrb_t
\quad \text{and}\quad  
\Qq_t - \Qqb_t 
	= - \Pp_t[t(\Dirp - \Dirb)]\Ppb_t 
		- \Qq_t[t(\Dirp - \Dirb)]\Qqb_t.
\end{equation}
Using the  hypothesis  that
$\Dirb - \Dirp = A_1 \conn + \divv A_2 + A_3$, 
\begin{equation}
\label{Eqn:PsiToPtQt}
\begin{aligned} 
&\norm{(\Qq_t - \Qqb_t)f}\\ 
	&\qquad\leq \norm{\Pp_t(tA_1\conn)\Ppb_t f } + \norm{\Pp_t(t\divv A_2)\Ppb_t f} + \norm{\Pp_t(tA_3)\Ppb_t f} \\
	&\qquad\qquad+\norm{\Qq_t(tA_1\conn)\Qqb_t f} + \norm{\Qq_t(t\divv A_2)\Qqb_t f} + \norm{\Qq_t(tA_3)\Qqb_t f}.
\end{aligned}
\end{equation}

\subsection{Reduction to quadratic estimates}

The goal of this subsection is to prove the 
following reduction of the main estimate
in Theorem \ref{Thm:Main} to the two 
quadratic estimates appearing the the hypothesis
of the following  proposition.
It is these two quadratic estimates that allow us to access 
real-variable harmonic analysis methods. The proofs
of these estimates are given in \S\ref{Sec:Harm1}
and \S\ref{Sec:Harm2}  respectively. 

\begin{proposition}
\label{Prop:FinalRed}
Suppose that 
\begin{align*}
&\int_{0}^1 \norm{\Qq_t A_1 \conn (\imath\iden + \Dir)^{-1} \Ppb_t u}^2\ \dtt 
	\leq C_1  \norm{A}_{\infty}^2 \norm{u}^2
\ \text{and} \ \\
&\int_{0}^1 \norm{t\Pp_t\divv A_2 \Ppb_t u}^2\ \dtt 
	\leq C_2 \norm{A}_{\infty}^2 \norm{u}^2
\end{align*}
for all $u \in \Lp{2}(\cV)$.
Then, for $\omega \in (0, \pi/2)$ and $\sigma \in (0, \infty)$,
whenever $f \in \Hol^\infty(\OSec{\omega,\sigma})$,
we obtain that 
$$\norm{f(\Dirp) - f(\Dirb)} \lesssim \norm{f}_\infty \norm{A}_\infty$$
where the implicit constant
depends  on $C_1,\ C_2$ and $\const(\cM,\cV,\Dirb,\Dirp)$.
\end{proposition}
 
First, we show that $f(\Dirb) \simeq f(\Dirp)$
can be reduced to a quadratic estimate
involving the difference of $\Qqb_t$ and $\Qq_t$. 
This is done via \eqref{Eqn:PsiToPtQt} and 
we estimate each of these terms using 
Proposition 4.5 and Proposition 4.7 in \cite{BMcR}.
Unlike in the situation of \cite{BMcR} where
the boundary was empty, we use the following
trace lemma to control the estimate on the boundary.
In what is to follow, 
$\Tr: \SobH{1}(\cV) \to \SobH{\frac{1}{2}}(\cW)$ is the boundary trace map. 

\begin{proposition}
\label{Prop:DivEst}
Let $\Uu_t$ be one of $\Rr_t$, $\Pp_t$ or $\Qq_t$
and $\Uub_t$ be one of $\Rrb_t$, $\Ppb_t$, $\Qqb_t$. 
 Then, 
$$ \sup_{t > 0} \norm{t \Uu_t \divv A_2 \Uub_t} \lesssim \norm{A_2}_\infty.$$  
\end{proposition}
\begin{proof}
Fix $u, v \in \Ck[c]{\infty}(\cV;\cB)$ and note that
$$\mh(\divv A_2 u, v) =  \mh(A_2 u, \conn v) + \divv W(u,v),$$ 
where $W(u,v) = (A_2)^j_{ik} u^i \delta_{jl} v^l\ dx^k$ 
inside an orthonormal frame, readily checked to be 
a well defined covectorfield. By Stokes' theorem,
$$\inprod{\divv A_2 u, v} - \inprod{A_2 u, \conn v} 
	= \int_{\Sigma} \mg(W(u,v)\rest{\Sigma}, \on)\ d\sigma.$$
By Cauchy-Schwartz, compactness of $\Sigma$
and smoothness of $\on$, we obtain that
$$\lmodulus{\int_{\Sigma} \mg(W(u,v)\rest{\Sigma}, \on)\ d\sigma} 
	\lesssim \norm{A_2}_\infty \norm{\Tr u} \norm{\Tr v}.$$

Next, note that whenever $\phi \in \dom(\Dirb)$
we have that $\phi \in \dom(\divv A_2 )$ and there exists a sequence $\phi_n \in \Ck[c]{\infty}(\cV; \B)$
such that $\phi_n \to \phi$ in $\dom(\Dirb)$ by the essential self-adjointness of $\Dirb$.
We prove that $\phi_n \to \phi$ in $\dom(\divv A_2)$. 
To prove this, note that $A_2: \Ck{\infty}(\cV) \to \Ck{0,1}(\cotanb\cM \tensor \cM)$
and fix a point $x \in \cM$, choose an orthonormal frame $\set{e_i}$ for $\cV$
and $\set{dx^i}$ for $\cotanb\cM$ with $\conn e_i = \conn dx^i = 0$ at $x$.
For $\psi \in \Ck{\infty}(\cV)$,
$A_2 \psi  = (A_2)_{i}^{jk} \psi^i\ dx^k \tensor e_j$, and 
$$ \divv A_2 \psi  
	= -\tr \conn \cbrac{ (A_2)_{i}^{jk}  \psi^i\ dx^k \tensor e_j}
 	= \sum_{k} (\partial_k (A_2)_{ik}^{j}) \psi^i + \sum_{k} (A_2)_{ik}^{j} \partial_k \psi^i)\  e_j.$$
Thus, $\modulus{\divv A_2 \psi}^2 \lesssim  \norm{\conn A_2}_\infty^2 \modulus{\psi}^2 
	+ \norm{A_2}_\infty^2 \modulus{\conn \psi}^2.$
Now, writing $\psi = \phi_n - \phi_m$, we obtain that
$$ \norm{\divv A_2(\phi_n - \phi_m)}^2 
	\lesssim \norm{\conn A_2}_\infty^2 \norm{\phi_n - \phi_m}^2 
	+ \norm{A_2}_\infty^2  \norm{\conn(\phi_n - \phi_m)}^2.$$ 
Since $\phi_n \in \dom(\Dirb)$, we have that  $\norm{\conn (\phi_n - \phi_m)} \lesssim \norm{\Dirb (\phi_n - \phi_m)} + \norm{\phi_n - \phi_m}$.
Thus, we have that $\phi_n \to \phi$ and $\divv A_2 \phi_n \to v$ and since $\divv A_2$ is closed as $A_2$ is bounded, we obtain
$\phi \in \dom(\divv A_2)$ and $v = \divv A_2 \phi$.

Now, let $u, v \in \Lp{2}(\cV)$.
Since we assume that $\Dir$ is essentially self-adjoint
on $\Ck[c]{\infty}(\cV;\cB)$, there exist sequences
$u_n, v_m \in \Ck[c]{\infty}(\cV; \cB)$ 
such that $u_n \to \Uub_t u$ and $v_m \to \Uu_t v$, 
with convergence in $\dom(\Dir)$, $\dom(\conn)$ and $\dom(\divv A_2)$
by what we have already established. Thus, 
\begin{align*}
\lmodulus{\inprod{t \Uu_t \divv A_2 \Uub_t u, v}} 
	&= \modulus{\lim_{m, n \to \infty}  \inprod{t \divv A_2 u_n, v_m}}  \\
	&\leq \lim_{m,n \to \infty} \modulus{\inprod{t A_2 u_n, \conn v_m}} 
		+ \lim_{m,n \to \infty}   \norm{A_2}_\infty   t \norm{\Tr u_n} \norm{\Tr v_m} \\
	&\lesssim \lim_{m,n \to \infty} \norm{A_2}_\infty \norm{u_n} (\norm{t \Dirp v_m} + t\norm{v_m}) \\
		&\qquad\qquad+   \norm{A_2}_\infty  t \norm{\Tr\Uub_t u} \norm{\Tr \Uu_t v} \\
	&\lesssim \norm{A_2}_\infty (\norm{u} + \sqrt{t} \norm{\Tr \Uub_t u}) \norm{v}, 
\end{align*}
where the last inequality follows from the standard boundary trace inequality.
on $\sqrt{t} \norm{\Tr \Uu_t v}$ and from the uniform
bounds on  $\norm{t \conn \Uu_t v} \lesssim \norm{t\Dirp \Uu_t v} + \norm{t \Uu_t v}$  and $t \norm{\Uu_tv }$.
 We obtain the conclusion by estimating $\norm{\Tr \Uub_t u}$ similarly.  
\end{proof}

As a consequence of this proposition  and \eqref{Eqn:PsiToPtQt},
we obtain  
$$\sup_{t \in (0,1]} \norm{\Uu_t - \Uub_t} \lesssim \norm{A}_\infty.$$  
Using this, arguing exactly as in \S4.2 in \cite{BMcR}, we
can reduce the required estimate in the conclusion of  Proposition \ref{Prop:FinalRed}
to proving a quadratic estimate: 
$$ \int_{0}^1 \norm{(\Qq_t - \Qqb_t)u}^2\ \dtt \lesssim \norm{A}_\infty^2 \norm{u}^2$$
for all $u \in \Lp{2}(\cV)$. 
 From \eqref{Eqn:PsiToPtQt},   
we obtain that 
\begin{equation}
\begin{aligned} 
&\cbrac{\int_{0}^1 \norm{(\Qq_t - \Qqb_t)u}^2\ \dtt}^{\frac{1}{2}}\\ 
	&\qquad\leq \cbrac{\int_0^1 \norm{\Pp_t tA_1\conn \Ppb_t u }^2\ \dtt}^{\frac{1}{2}}
		+ \cbrac{\int_0^1 \norm{\Pp_t t\divv A_2 \Ppb_t u}^2\ \dtt}^{\frac{1}{2}} \\
		&\qquad\qquad\qquad\qquad\qquad+ \cbrac{\int_0^1 \norm{\Pp_t tA_3 \Ppb_t u}^2\ \dtt}^{\frac{1}{2}} \\
		&\qquad\qquad+
		\cbrac{\int_0^1 \norm{\Qq_t tA_1\conn \Qqb_t u}^2\ \dtt}^{\frac{1}{2}}
\label{Eqn:PsiToPtQt2}
		+ \cbrac{\int_0^1 \norm{\Qq_t t\divv A_2 \Qqb_t u}^2\ \dtt}^{\frac{1}{2}}\\ 
		&\qquad\qquad\qquad\qquad\qquad+ \cbrac{\norm{\Qq_t tA_3 \Qqb_t u}^2\ \dtt}^{\frac{1}{2}}. 
\end{aligned}
\end{equation}

Estimating as in Proposition 4.7 in \cite{BMcR}, we bound
the first, third and sixth term by $\norm{A}_\infty^2 \norm{f}^2.$
The second and forth terms are controlled by the 
hypothesis of Proposition \ref{Prop:FinalRed}.
The only term that remains to be bounded is the penultimate term in this expression
for which the estimate in Proposition 4.7 in \cite{BMcR} does not work.  
The way in which we estimate this term requires a slight excursion into interpolation theory.

Let $\SobH{1}(\cV)$ denote the first-order Sobolev space on $\cV$
and define 
$$\SobH{s}(\cV) = \cinterpol{\Lp{2}(\cV), \SobH{1}(\cV)}{\theta = s},$$
for $s \in [0,1]$ where $\cinterpol{\cdot, \cdot}{\theta}$ represents complex interpolation.
Also, let 
$$\SobH[0]{s}(\cV) = \close{\Ck[cc]{\infty}(\cV)}^{\norm{\mdot}_{\SobH{s}}}, 
\ \SobH{-s}(\cV) = \adj{\SobH[0]{s}(\cV)},\quad\text{and}\quad
\SobH[00]{s}(\cV) = \cinterpol{\Lp{2}(\cV), \SobH[0]{1}(\cV)}{\theta = s}.$$

In order to gain an explicit expression for the norms in these
interpolation scales, we connect these spaces to domains 
of operators. 
Let $\conn[N] = \close{\conn[2]}$ and $\conn[D] = \close{\conn[0]}$,
where $\conn[2]: \Ck{\infty}\intersect \Lp{2}(\cV) \to \Ck{\infty}\intersect \Lp{2}(\cotanb\cM \tensor \cV)$
and $\conn[0]: \Ck[cc]{\infty}(\cV) \to \Ck[cc]{\infty}(\cotanb\cM \tensor \cV)$.
The subscripts ``$N$'' and ``$D$'' are chosen for Neumann and Dirichlet respectively
since $\SobH{1}(\cV) = \dom(\conn[N]) = \dom(\sqrt{\Lap_N})$ and 
$\SobH[0]{1} = \dom(\conn[D]) = \dom(\sqrt{\Lap_D})$, where
$\Lap_N = \adj{\conn[N]}\conn[N]$ and $\Lap_D = \adj{\conn[D]}\conn[D]$.
Moreover, $\norm{\mdot}_{\SobH{1}} \simeq \norm{(\iden + \sqrt{\Lap_N})\mdot}$
and $\norm{\mdot}_{\SobH[0]{1}} \simeq \norm{(\iden + \sqrt{\Lap_D})\mdot}$. 
 
Consequently,by Theorem 6.6.9 in \cite{Haase},  we have that:
\begin{align*}
\SobH{s}(\cV) &= \cinterpol{\Lp{2}(\cV), \SobH{1}(\cV)}{\theta = s} = \dom((\iden + \sqrt{\Lap_N})^s),\\
\SobH[00]{s}(\cV) &=\cinterpol{\Lp{2}(\cV), \SobH[0]{1}(\cV)}{\theta = s} =  \dom((\iden + \sqrt{\Lap_D})^s), 
\end{align*}
and in particular for $s \in [0,1]$,
$$\norm{\cdot}_{\SobH{s}} \simeq \norm{ (\iden + \sqrt{\Lap_N})^s \cdot}
\quad \text{and}\quad
\norm{\cdot}_{\SobH{-s}} \simeq \norm{ (\iden + \sqrt{\Lap_N})^{-s} \cdot}.$$
Since the identity map embeds $\SobH[00]{1}(\cV) \embed \SobH{1}(\cV)$
and $\SobH[00]{0}(\cV) \embed \SobH{0}(\cV)$, we have 
by interpolation that 
$$\dom((\iden + \sqrt{\Lap_D})^s) = \SobH[00]{s}(\cV)  \embed \SobH{s}(\cV) = \dom((\iden + \sqrt{\Lap_N})^s)$$
for $s \in (0,1)$. 
Similarly, since $\dom(\Dir) = \dom(\modulus{\Dir})$, where $\modulus{\Dir} = \sqrt{\Dir^2}$
and $\norm{(\iden + \modulus{\Dir})u}  \simeq \norm{u} + \norm{\Dir u}$, by
the same Theorem 6.6.9 in \cite{Haase}, 
$$ \cinterpol{\Lp{2}(\cV), \dom(\Dir)}{\theta = s} = \dom(\modulus{\Dir}^s) = \dom( (\iden + \modulus{\Dir})^s).$$

 The following key result is well known in the case of functions on the 
upper half space and smooth Euclidean domains by the work 
of Bergh and Löfström in \cite{BL} or Triebel in \cite{Triebel}.
The following is a vector bundle version which, to our knowledge,
does not seem to have been treated previously in the literature. 
\begin{lemma}
\label{Lem:SmallSob}
The equality $\SobH{s}(\cV) = \SobH[0]{s}(\cV) = \SobH[00]{s}(\cV)$ holds whenever
 $0 \leq s < 1/2$.
\end{lemma}
\begin{proof}
Now let $U_0 = \cM \setminus Z$, where $Z$ is a smooth precompact open
neighbourhood of $\Sigma = \bnd\cM$ and $(\phi_j, \psi_j, U_j)$
trivialisations $\psi_j$ inside charts $\phi_j:U_j \to \R^{n}_+$
for $j = 1, \dots, M$, so that $M = \union_{j=0}^M U_j$.
Let $\set{\eta_j}$ be a smooth partition of unity
subordinate to $\set{U_j}$. We can choose $\eta_j$ such that $\modulus{\conn \eta_j} \leq C$
for some $C > 0$.

Define:
\begin{align*}  
&B_0 = \Lp{2}(\cV), &&A_0 = \Lp{2}(\cV) \oplus \Lp{2}(\R^{n}_+, \C^N)^M \\
&B_1 = \SobH{1}(\cV), &&A_1 = \SobH[0]{1}(\cV) \oplus \SobH{1}(\R^n_+,\C^N)^M \\
&B_1^0 = \SobH[0]{1}(\cV), &&A_1^0 = \SobH[0]{1}(\cV) \oplus \SobH[0]{1}(\R^n_+, \C^N)^M.
\end{align*}
Now, define 
$S: B_0 \to A_0$ by
$$Su = (\eta_0, \psi_1(\eta_1 u)\comp \phi_1^{-1}, \dots, \psi_M(\eta_M u) \comp \psi^{-1}_M),$$
with $j$-th coordinate map extended to $0$ outside of the support of $\eta_j$, 
and note $S$ is an injection. Moreover,
it is also a map $B_1 \mapsto A_1$ and $B_1^0 \mapsto A_1^0$.
Also, define $R: A_0 \to B_0$ by 
$$R(u_0, u_1, \dots, u_M) = u_0 + \eta_1 \psi_1^{-1}( u_1 \comp \phi_1 ) + \dots +  \eta_M \psi_M^{-1}( u_M \comp \phi_M).$$
It is also easy to see that this is a map $A_1 \mapsto B_1$ and $A_1^0 \mapsto B_1^0$. 

Now, note that $RS = \iden$ on $\bddlf(B_j, B_j)$ for $j=0,1$  and  $\bddlf(B_1^0, B_1^0)$.
That is, $R$ is a \emph{retraction} and $S$ is a \emph{coretraction}
associated to $R$. By Theorem ($\ast$) in \S1.2.4 of  \cite{Triebel}
we get that $S$ is an isomorphic mapping from $\SobH{s}(\cV) \xrightarrow{\cong} W$
 for $s \in (0,1)$  
where $W$ is a  closed   subspace of $\SobH[00]{s}(\cV) \oplus \SobH{s}(\R^n_+, \C^N)^M$.
Similarly, we have that  $\SobH[00]{s}(\cV) \xrightarrow{\cong} W_0$  
with $W_0$ is a closed subspace of  $\SobH[00]{s}(\cV) \oplus \SobH[00]{s}(\R^n_+, \C^N)^M$. 
The subspace $W$ is the range of $SR$ restricted to  $\SobH[00]{s}(\cV) \oplus \SobH{s}(\R^n_+,\C^N)^M$ 
and similarly $W_0$ is the range of $SR$ restricted to  $\SobH[00]{s}(\cV) \oplus \SobH[00]{s}(\R^n_+,\C^N)^M$.  
 But by Theorems 11.1 and 11.2  in \cite{BL}, 
we obtain  $\SobH[0]{s}(\R^{n}_+,\C^N) = \SobH[00]{s}(\R^n_+, \C^N) = \SobH{s}(\R^n_+, \C^N)$   for  $0 \leq s < 1/2$,  
and therefore, $W_0 = W$ for $0 \leq s < 1/2$.  This shows
that $\SobH{s}(\cV) = \SobH[00]{s}(\cV)$ for $0 \leq s < 1/2$.

To finish off the proof, note that
$\norm{(\iden + \sqrt{\Lap_N})u} \lesssim \norm{(\iden + \sqrt{\Lap_D})u}$
so through interpolation we get
$\norm{(\iden + \sqrt{\Lap_N})^s u} \lesssim \norm{(\iden + \sqrt{\Lap_D})^s u}$.
Since $\Ck[cc]{\infty}(\cV)$ is dense in $\SobH[00]{s}(\cV) = \dom((\iden + \sqrt{\Lap_D})^s)$,
we have that $\SobH[00]{s}(\cV) \embed \SobH[0]{s}(\cV)$.
But we have $\SobH[0]{s}(\cV) \embed \SobH{s}(\cV)$ and 
since we have already proved $\SobH{s}(\cV) = \SobH[00]{s}(\cV)$ for $0 \leq s < 1/2$,
we obtain the conclusion. 
\end{proof} 

With the aid of this lemma, we obtain the following.
\begin{proposition}
\label{Prop:QtQest} 
The quadratic estimate
$$ \int_0^1 \norm{\Qq_t t\divv A_2 \Qqb_t f}^2\ \dtt \lesssim \norm{f}^2$$
holds for $f \in \Lp{2}(\cV)$.
\end{proposition}
\begin{proof}
Fix $u \in \Lp{2}(\cV)$ and estimate
$$\inprod{\Qq_t t \divv A_2 \Qqb_t f, u} = - \inprod{A_2 \Qqb_t f, t \conn \Qq_t u} +  t\inprod{A_2 \Tr \Qqb_t f, \Tr \Qq_t u}_{\Lp{2}(\cW)}.$$ 
It is easy to see that
$$\modulus{\inprod{A_2 \Qqb_t f, t \conn \Qq_t u}} \lesssim \norm{A_2}_\infty \norm{u} \norm{\Qqb_t f},$$
so it remains to consider the  boundary term.
Note that
$$ \modulus{ t\inprod{A_2 \Tr \Qqb_t f, \Tr \Qq_t u}_{\Lp{2}(\Sigma)}} 
	\lesssim \norm{A_2}_\infty t \norm{\Tr \Qqb_t f}_{\Lp{2}(\cW)} \norm{\Tr \Qq_t u}_{\Lp{2}(\cW)}.$$ 
By the standard boundary trace inequality, we obtain that
$\sqrt{t} \norm{\Tr \Qq_t u}_{\Lp{2}(\cW)} \lesssim \norm{u}$.

To bound $\Qqb_t f$, let $\vec{N}$ be an extension of the  normal vectorfield $\on$ 
on a compact neighbourhood around $\Sigma$.
Then, 
\begin{align*}
 t\norm{\Tr \Qqb_t f}^2_{\Lp{2}(\cW)}   &= t \int_{\cM} \divv (\modulus{\Qqb_t f}^2 \vec{N})\ d\mu \\
	&\lesssim t \int_{\cM} \re \mg(\conn[\vec{N}] \Qqb_t f, \Qqb_t f)\ d\mu + t \norm{\Qqb_t f}^2 \\ 
	&\lesssim t \modulus{\inprod{\conn[\vec{N}] \Qqb_t f, \Qqb_t f}} + t \norm{\Qqb_t f}^2.
\end{align*}
  On fixing $0 < s < 1/2$,  we note that
\begin{equation}
\label{Eq:IntReq}
\modulus{\inprod{\conn[\vec{N}]\Qqb_t f, \Qqb_t f}} \lesssim 
	\norm{\conn[\vec{N}] \Qqb_t f}_{\SobH{-s}} \norm{\Qqb_t f}_{\SobH{s}},
\end{equation}
Now, note that $\conn[\vec{N}]:\SobH{1}(\cV) \to \Lp{2}(\cV)$ and on defining
$(\conn[\vec{N}]u)(v) = - \inprod{u, \conn[\vec{N}]v}$ for $v \in \Ck[c]{\infty}(\cV)$,
we obtain that $\conn[\vec{N}]: \Lp{2}(\cV) \to \adj{\SobH[0]{1}(\cV)}  = \SobH{-1}(\cV)$
 boundedly. 
By interpolation, we obtain 
that 
$\conn[\vec{N}]: \cinterpol{\SobH{1}(\cV), \Lp{2}(\cV)}{\theta = s} 
	\to \cinterpol{\Lp{2}(\cV), \SobH{-1}(\cV)}{\theta = s}$ boundedly.
Note, however, that 
$$\cinterpol{\SobH{1}(\cV), \Lp{2}(\cV)}{\theta = s} = \cinterpol{\Lp{2}(\cV), \SobH{1}(\cV)}{\theta = 1 - s} = \SobH{1-s}(\cV),$$
and  that  
$$
\cinterpol{\Lp{2}(\cV), \SobH{-1}(\cV)}{\theta = s}
=  \adj{(\cinterpol{\Lp{2}(\cV), \SobH[0]{1}(\cV)}{\theta = s})} 
= \adj{\SobH[00]{s}(\cV)}
= \adj{\SobH[0]{s}(\cV)} 
= \SobH{-s}(\cV),$$   
where we have used that $\Lp{2}(\cV)$ is reflexive 
and Corollary 4.5.2 in \cite{BL}
in the first equality and that 
$s < 1/2$ and Lemma \ref{Lem:SmallSob}
in the penultimate equality. 
On combining these facts, we obtain that 
\begin{equation*} 
\norm{\conn[\vec{N}]\Qqb_t f}_{\SobH{-s}} \lesssim  \norm{\Qqb_t f}_{\SobH{1 -s}}.
\end{equation*}
Moreover, since $\dom(\modulus{\Dir}) \embed \SobH{1}(\cV)$ and $\dom(\modulus{\Dir}^0) = \Lp{2}(\cV) \embed \SobH{0}(\cV) = \Lp{2}(\cV)$,
we have   $\dom(\modulus{\Dir}^q) \embed \SobH{q}(\cV)$  for $q \in [0,1]$ by interpolation
and hence, 
$$t^q \norm{\Qqb_t f}_{\SobH{q}} \lesssim \norm{t^q(\iden + \modulus{\Dirb}^q) \Qqb_t f} 
	\leq \norm{\psi_{q}(t\Dirb) f} + \norm{\Qqb_t f},$$
where $\psi_q(\zeta) = \zeta \modulus{\zeta}^q(1 + \zeta^2)^{-1}.$
Thus, 
\begin{multline*}
t \modulus{\inprod{\conn[\vec{N}]\Qqb_t f, \Qqb_t f}} \\ \lesssim   
	(t^{1 - s} \norm{\Qqb_t f}_{\SobH{1-s}}) (t^s \norm{\Qqb_t f}_{\SobH{s}}) 
	\lesssim \norm{\psi_{1 - s}(t\Dirb)f}^2 + \norm{\psi_{s}(t\Dirb) f}^2 + \norm{\Qqb_t f}^2, 
\end{multline*}  
and therefore, 
$$ 
t\norm{\Tr \Qqb_t f}_{\Lp{2}(\cW)}  \lesssim 
	 \norm{\psi_{1 - s}(t\Dirb)f}^2 + \norm{\psi_{s}(t\Dirb) f}^2 + (1 + t)\norm{\Qqb_t f}^2.$$
Noting that 
$$\int_{0}^1 \norm{\psi_{q}(t\Dirb)f}^2\ \frac{dt}{t} \leq C_q \norm{f}^2$$
for $q \in [0,1)$ completes the proof. 
\end{proof} 

\begin{remark}
The equation \eqref{Eq:IntReq} demonstrates the necessity of the interpolation methods since we can only conclude the desired quadratic estimates provided a derivative of order strictly less than $1$ is applied to $\Qqb_t f$.
\end{remark}

\subsection{Harmonic analysis I}
\label{Sec:Harm1} 

In this subsection, on drawing from the 
estimates in \S5 in \cite{BMcR}, we
demonstrate how to handle the first
quadratic estimate term
$$\int_{0}^1 \norm{\Qq_t A_1 \conn (\imath\iden + \Dir)^{-1} \Ppb_t f}^2\ \dtt 
	\lesssim   \norm{A}_{\infty}^2 \norm{f}^2 $$
appearing in the hypothesis of
Proposition \ref{Prop:FinalRed}.
In order to avoid repetition, we encourage
the reader to keep a copy of \cite{BMcR} 
handy to navigate through the remainder of this paper.

The following is an itemisation of the notation 
that we will require from \S5 of \cite{BMcR}:  
\begin{itemise} 
\item  Dyadic cubes  $\set{\Q[\alpha]^k \subset \cM: \alpha \in I_k,\ k \in \Na}$,
with centres $z_\alpha^k \in \Q[\alpha]^k$, where $\union_{k} \Q[\alpha]^k$ cover
$\cM$ almost everywhere, and when $\beta > \alpha$, $\Q[\alpha]^k \intersect \Q[\beta]^l = \emptyset$
or $\Q[\alpha]^k \subset \Q[\beta]^l$.
The cubes are of a fixed ``length'' $\delta \in (0,1)$, 
and a $\delta^j$ cube contains an $a_0 \delta^j$ ball
and has diameter at most $C_1 \delta^j$.
The length of a cube $\Q$ is denoted $\len(\Q)$. 
The constant $\eta > 0$ is an exponent that measures smallness of 
the volume toward the edge of a cube with constant $C_2 > 0$. 
See Theorem 5.1 in \cite{BMcR}.

\item The scale is defined as $\scale = \delta^\jscale$ 
where  $C_1 \delta^{\jscale} \leq \brad/5$, 
with $\brad = \max\set{\brad_{\cotanb\cM}, \brad_{\cV}}$,
the maximum of the GBG radii of $\cotanb\cM$ 
and $\cV$.

\item The collection of dyadic cubes $\DyQ^j$, 
$\DyQ = \union_{j \geq \jscale} \DyQ^j$, and $\DyQ_t$
for $t \leq \scale$.

\item The unique ancestor $\ancester{\Q} \in \DyQ^J$ for a dyadic cube $\Q$, 
 the set of \emph{GBG coordinates} $\sC$, which for a cube $\Q \in \DyQ^j$ 
is the GBG trivialisation pertaining to the unique GBG ball containing the cube in $\DyQ^J$ containing $\Q$,
and \emph{dyadic GBG coordinates} $\sC_\jscale$ which is the restriction of 
this GBG ball to the cube which contains it.

\item The \emph{cube integral}
	$\Ball(x_{\ancester{\Q}},\brad) \times \DyQ \ni (x, \Q) \mapsto (\int_{\Q} \mdot )(x)$
	defined on $\Lp[loc]{1}(\cV)$
	by 
	$$\cbrac{\int_{\Q} u}(x)= \cbrac{\int_{\Q} u^i(y)\ d\mu(y)} e_i(x)$$
	where $e_i$ is the GBG coordinates of $\Q$, and \emph{cube average}  $u_{\Q} = \fint_{\Q} u$ inside 
	the GBG coordinate ball of $\Q$ and $0$ outside it. 

\item For $t > 0$, the \emph{dyadic averaging operator}
$\Av_t: \Lp[loc]{1}(\cV) \to \Lp[loc]{1}(\cV)$
	given by $\Av_t(x) = (\fint_{\Q} u)(x)$ where $x \ni \Q$.

\item For a $w = w^i e^{\C^N}_i \in \C^N$, the locally constant
	extension inside the GBG coordinates of $\Q$ 
	are given by $\omega^c(x) = w^i e_i(x)$ and 
	zero outside of this coordinate ball.

\item Given a $t$-uniformly bounded family of operators $\QQ_t$,
	define the principal part 
	$\Pri^\QQ_t(x): \C^N \cong \cV_x \to \cV_x$ of $\QQ_t$ by  
	by $\Pri^\QQ_t(x) w = (\QQ_t \omega^c)(x)$.
\end{itemise} 

The following is a key lemma that
is necessary in order to adapt the arguments 
of  \S5 of \cite{BMcR}  to our manifold with boundary.
It allows us to ensure that we can use a cutoff that restricts
the estimates \emph{away} from the boundary.

\begin{lemma}
\label{Lem:CanLem}
There exist constants $k_0, \tilde{\eta}, \tilde{C}_3 > 0$ such that
for all cubes $\Q \in \DyQ^{k}$ with $k >  k_0$ and $\close{\Q} \intersect \Sigma \neq \emptyset$,
we have 
$$\mu\set{x \in \Q:\met(x, \Sigma) \leq s \len(\Q)} \leq \tilde{C}_3 s^{\tilde{\eta}} \mu(\Q).$$
In particular, for every $\Q \in \DyQ^k$ with $k > k_0$, 
$$\mu\set{x \in \Q: \met(x, \cM\setminus( \Q \setminus \Sigma)) \leq s \len(\Q)}
	\leq \tilde{C}_3 s^{\tilde{\eta}} \mu(\Q).$$
The constants $\tilde{\eta}$ 
and $\tilde{C}_3$ depends on $\eta, a_0$ and $C_1$
from Theorem 5.1 in \cite{BMcR}.
\end{lemma}
\begin{proof}
Let $Z = \set{x \in \cM: \met(x, \Sigma) \leq \epsilon}$ with $\epsilon < 1$ chosen sufficiently small so that $Z$ is a smooth compact submanifold  of $\cM$ with smooth boundary $\Sigma$.
Let $\tilde{Z}$ be the smooth compact manifold without boundary obtained
by taking two copies of $Z$ and identifying the boundaries, 
and extending the metric appropriately. This metric is $\Ck{0}$
and there exists a smooth $\Ck{\infty}$ metric $G$-close to $\mg$
for some $G \geq 1$. Consequently, without loss
of generality, we assume that the metric extension is smooth.
Let $k_{\Sigma} = \inj(\tilde{Z}) > 0$.

By the compactness of $\tilde{Z}$, 
we use Theorem 1.2 in \cite{Hebey} to obtain 
$C_\Sigma \geq 1$ such that for each $x \in \tilde{Z}$,
$(\psi_x, B(\frac{1}{2}k_\Sigma, x))$ is a coordinate
chart with
$$ C_\Sigma^{-1} \modulus{u}_{\pullb{\psi}_x\delta(y)}
	\leq \modulus{u}_{\mg(y)} \leq C_\Sigma \modulus{u}_{\pullb{\psi}_x\delta(y)},$$
for each $y \in B(\frac{1}{2}k_\Sigma, x)$, and where $\delta$ is the Euclidean metric
in that chart. In particular, since $Z \subset \tilde{Z}$
and the topology of $Z$ is the subspace topology inherited
from $\tilde{Z}$, we get that this holds for balls $B(x,r)$
in $Z$ as well.
From this, inside $(\psi_x, B(\frac{1}{2}k_\Sigma,x))$,
on letting $\met^\ast(x,y) = |\psi_x(x) - \psi_y(y)|$
and $\Leb[\ast] = \pullb{\psi_x}\Leb$,
\begin{equation}
\label{Eq:CanEst}
C_\Sigma^{-1} \met^\ast(x,y) \leq \met(x,y) \leq C_\Sigma \met^\ast(x,y)\quad \text{and}\quad  
C_\Sigma^{-\frac{n}{2}}\ d\Leb[\ast] \leq d\mu \leq C_{\Sigma}^{\frac{n}{2}}\ d\Leb[\ast].
\end{equation}
Now, fix $k_0 > 0$ such that so that $C_1 \delta^{k_0} < \frac{1}{10}k_\Sigma$.
Then, for all $k > k_0$, whenever $\Q \in \DyQ^k$, we have that 
$\Q \subset B(x_{\Q}, \frac{1}{2}k_\Sigma)$, which corresponds
to a coordinate system with control on the metric and measure
as we have describe before. 

Fix such a cube $\Q \in \DyQ^k$  and 
define $\Q[\Sigma, s] = \set{x \in \Q: \met(x,\Sigma) \leq s \len(\Q)}$
and note that on using \eqref{Eq:CanEst},  
$$ \psi_{\Q}(\Q[\Sigma,s]) \subset E_{\Sigma,s} 
	= \set{x \in \psi_{\Q}(\Q): \met_{\R^n}(x, \R^{n-1} 
	\intersect \close{\psi_{\Q}(\Q)} \leq C_\Sigma s \delta^k}.$$
Similarly, we have that
$\psi_{\Q}(B(x_{\Q},C_1 \delta^k)) \subset B_{\R^n}(\bar{x}_{\Q}, C_\Sigma C_1 \delta^k)
\subset \mathrm{Box}_{\R^n}(\bar{x}_{\Q}, C_\Sigma C_1 \delta^k)$
where $\bar{x}_{\Q} = \psi_{\Q}(x_{\Q})$
and $\mathrm{Box}_{\R^n}(x,l)$ is a Euclidean 
box centred at $x$ of length $l$. Then,
\begin{multline*} 
\Leb(E_{\Sigma,s}) 
	\leq \Leb[n-1](\R^n \intersect \mathrm{Box}_{\R^n}(\bar{x}_{\Q}, C_{\Sigma}C_1\delta^k)) 
		\times C_\Sigma s \delta^k \\
	\leq (C_\Sigma C_1 \delta^k)^{n-1} \times C_{\Sigma} s \delta^k
	= C_\Sigma^n C_1^{n-1} \delta^{nk} s.
\end{multline*}

Similarly, we have that $\psi_{\Q}(B(x_{\Q}, a_0 \delta^k)) \supset B_{\R^n}(\bar{x}_{\Q}, C_{\Sigma}^{-1} a_0 \delta^k)$, 
and
\begin{multline*}
\frac{\mu(\Q[\Sigma, s])}{\mu(\Q)} 
	\leq \frac{\mu(\Q[\Sigma, s])}{\mu(B(x_{\Q}, a_0 \delta^k)}
	\leq \frac{C_{\Sigma}^{\frac{n}{2}} \Leb(E_{\Sigma,s})}
		{C_{\Sigma}^{-\frac{n}{2}}\Leb(B_{\R^n}(\bar{x}_{\Q}, C_{\Sigma}^{-1} a_0 \delta^k))} \\
	\leq C_{\Sigma}^n \frac{C_\Sigma^n C_1^{n-1} \delta^{nk} s}
		{\omega_n (C_{\Sigma}^{-1}  a_0 \delta^k)^n}
	= \frac{C_{\Sigma}^{3n} C_1^{n-1}}{\omega_n a_0^n} s, 
\end{multline*}
where the first estimate follows from Theorem 5.1 (v) in \cite{BMcR},
the second estimate from our previous calculation
combined with \eqref{Eq:CanEst}, and where $\omega_n$ is the 
volume of the ball of unit radius in $\R^n$.

Set 
$\tilde{\eta} = \max\set{1, \eta}$ and
$\tilde{C}_3 = \max\set{C_3,  \frac{C_{\Sigma}^{3n} C_1^{n-1}}{\omega_n a_0^n}},$
and noting
$$\set{x \in \Q: \met(x, \cM\setminus( \Q \setminus \Sigma)) \leq s \len(\Q)}
	= \set{x \in \Q: \met(x, \cM\setminus\Q) \leq s \len(\Q)} \union \Q[\Sigma,s],$$
completes the proof.
\end{proof}

\begin{proposition}
\label{Prop:HarmAnal1} 
The quadratic estimate 
$$\int_{0}^1 \norm{\Qq_t A_1 \conn (\imath\iden + \Dir)^{-1} \Ppb_t u}^2\ \dtt 
	\lesssim \norm{A}_{\infty}^2 \norm{u}^2$$
holds for all $u \in \Lp{2}(\cV)$, with the implicit
constant depending on $\const(\cM,\cV, \Dirb, \Dirp)$.
\end{proposition} 
\begin{proof}
We split the estimate as follows: 
\begin{align*}
\int_0^1 \norm{\Qq_t A_1 \conn (\imath \iden + \Dir)^{-1} \Ppb_t u}^2\ \dtt
	&\lesssim  \int_0^1 \norm{(\Qq_t - \Pri_t\Av_t) \A_1 \conn (\imath \iden + \Dir)^{-1} \Ppb_t u}^2\ \dtt \\
		&\quad+ \int_0^1 \norm{\Pri_t\Av_t\A_1 \conn (\imath \iden + \Dir)^{-1}(\iden - \Ppb_t) u}^2\ \dtt  \\
		&\quad+\int_0^1 \norm{\Pri_t \Av_t \A_1 \conn (\imath \iden + \Dir)^{-1} u}^2\ \dtt.
\end{align*}

Now, we note that the off-diagonal decay given in Lemma 5.9
in \cite{BMcR} is valid for our operator $\Qq_t A_1$  due to the
local boundary conditions encoded in assumption \ref{Hyp:Comm}.  
Thus, we can apply  Propositions 5.4, Lemma 5.8 and Proposition 5.12 
in \cite{BMcR} to estimate the terms appearing in this
decomposition. We give a brief description of how this is done.

The first term is estimated by using an
argument similar to the proof of Proposition 5.4 and Theorem 2.4 in \cite{BMcR}, 
with $\cW = \cotanb\cM \tensor \cV$. It suffices to note that
since $\norm{u}_{\Dirb} \simeq \norm{u}_{\SobH{1}}$ for $u \in \dom(\Dirb)$,
this argument can be run in verbatim. It simply remains
to prove $\norm{\conn S u} \lesssim \norm{u}_{\SobH{1}}$ for $S =  \conn(\imath \iden + \Dir)^{-1}$.
This argument is included in the proof of Theorem 2.4 on noting that the argument
runs in verbatim due to  assumption \ref{Hyp:Weitz}. 
 
For the middle term in the estimate, we use the argument
in proving Proposition 5.10 in \cite{BMcR}. This argument
is straightforward from establishing the \emph{cancellation lemma}, Lemma 5.8 in \cite{BMcR}. 
To prove this lemma, we note that for each dyadic cube $\Q$, 
and for each $u \in \dom(\Dirb)$ with $\spt u \subset \Q \intersect \interior{\cM}$,
we have that
$$ \left|{\int_{\Q} \Dirb u\ d\mu}\right| \lesssim \mu(Q)^{\frac{1}{2}} \norm{u}
\quad\text{and}\quad
\left|{\int_{\Q} \conn u\ d\mu}\right| \lesssim \mu(Q)^{\frac{1}{2}} \norm{u}, $$
where the implicit constants depends on $\const(\cM,\cV,\Dirb,\Dirp)$. 
On coupling these estimates with Lemma \ref{Lem:CanLem}, we obtain 
the statement of Lemma 5.8 in \cite{BMcR} in our present context.

The last term is obtained by a straightforward application of Proposition 5.12 in \cite{BMcR}.
\end{proof}

\subsection{Harmonic analysis II}
\label{Sec:Harm2}

In this subsection, we prove the remaining estimate 
$$\int_{0}^1 \norm{t\Pp_t\divv A_2 \Ppb_t u}^2\ \dtt 
	\lesssim \norm{A}_{\infty}^2 \norm{u}^2$$
for all $u \in \Lp{2}(\cV)$.
It is in the proof of this estimate where the main novelty of the harmonic analysis
in this paper can be found.
A key difficulty here is that the off-diagonal decay  - and even $\Lp{2}$-boundedness -  of $t \Pp_t\divv A_2$, 
which holds when $\cM$ has no boundary,  is not 
valid due to the fact that $A_2$  does not preserve  boundary conditions.
Despite this obstacle, on considering 
the operator $t\Pp_t \divv A_2 \Ppb_t$ instead as a whole,  we are able to prove 
the required quadratic estimate.
Our approach here is motivated by a similar argument in \cite{AAH}
by Auscher, Axelsson (\Rosen) and Hofmann. 

For the remainder of this subsection, 
let 
$$\Theta_t =t\Pp_t\divv A_2 \Ppb_t$$
and let $\Pri_t$ denote the principal part 
of $\Theta_t$ we recall is
$\Pri_t^\Theta(x) w = (\Theta_t \omega^c)(x)$, 
where $\omega^c$ is the constant section
related to $w \in \cV_x \cong \C^N$.

\begin{lemma}
\label{Lem:ThetaOD}
The operators  $\Theta_t$ are uniformly bounded in $t > 0$ and
have the off-diagonal decay estimate: 
there exists $C_{\Theta} > 0$ such that, for each $M > 0$, there exists a constant $C_{\Delta, M} > 0$
with
\begin{align*}
\norm{\ch{E} \Theta_t(\ch{F}u)}_{\Lp{2}(\cV)}
	\leq  C_{\Delta,M} \norm{A}_{\infty} \maxx{\frac{\met(E,F)}{t}}^{-M} 
		\exp\cbrac{-C_{\Theta} \frac{\met(E,F)}{t}} \norm{\ch{F}u}_{\Lp{2}(\cV)},
\end{align*}
for every Borel set $E,\ F \subset \cM$, $u \in \Lp{2}(\cV)$, and where
$\maxx{a} = \max\set{1,a}$.
\end{lemma}
\begin{proof}  
Uniform bounds for $\Theta_t$ were proved in Proposition \ref{Prop:DivEst}. 
Building on this, we prove the off-diagonal estimates in the
conclusion by reduction to corresponding such 
estimates for the resolvents $\Rrb_t$ and $\Rr_t$, 
which are immediate by replicating the argument of Lemma 5.3 in \cite{CMcM}
in light of \ref{Hyp:Comm}.

Given $E, F \subset \cM$ Borel with $\met(E,F) > 0$, 
pick $\eta \in \Ck{\infty}(\cM)$ such that 
$\eta(x) = 1$ when $\met(x,E) < 1/3 \met(E,F)$
and $\eta(x) = 0$ when $\met(x,F) < 1/3 \met(E,F)$ 
so that $\norm{\conn \eta}_\infty \lesssim 1/\met(E,F)$.
It suffices to prove the required estimates
for $\Rr_t t \divv A_2 \Rrb_t$ since by replacing
$t$ by $-t$ in the estimates below and noting $\Ppb_t = (\Rrb_t + \Rrb_{-t})/2$
and similarly $\Pp_t = (\Rr_t + \Rr_{-t})/2$ yields the bound 
for $\Theta_t$.
Now, note that
$$ 
\norm{ \ind{E} \Rr_t t \divv A_2 \Rrb_t (\ind{F} u)}
	= \norm{ \ind{E} \comm{ \eta, \Rr_t t \divv A_2 \Rrb_t} \ind{F} u}$$
and 
\begin{multline*}
\comm{ \eta, \Rr_t t \divv A_2 \Rrb_t} \\
	= - \Rr_t \comm{\eta, \imath t \Dirp} \Rr_t t \divv A_2 \Rrb_t
	+ \Rr_t \comm{\eta, t \divv} A_2 \Rrb_t
	- (\Rr_t t \divv A_2 \Rrb_t) \comm{\eta, \imath t \Dirb} \Rr_t.$$
\end{multline*}

Since $\comm{\eta, \Dirp},\ \comm{\eta, \divv},\ \comm{\eta, \Dirb}$ 
are multiplication operators whose $\Lp{\infty}$ norm is
bounded by $\norm{\conn \eta}_\infty$  and supported on 
$$G = \set{x \in \cM: \met(x,E) \geq \frac{1}{3} \met(E,F)\ \text{and}\ \met(x,F) \geq \frac{1}{3} \met(E,F)},$$
we obtain the conclusion from off-diagonal estimates
for $\Rr_t: \Lp{2}(G;\cV) \to \Lp{2}(E; \cV)$
and $\Rrb_t: \Lp{2}(F; \cV) \to \Lp{2}(G;\cV)$,
and from uniform bounds on $\Rr_t t \divv A_2 \Rrb_t$
from Proposition \ref{Prop:DivEst}.
\end{proof}

Next, we split the required estimate in the following way: 
\begin{equation}
\label{Eq:ThetaSplit}
\begin{aligned}
\int_{0}^1 \norm{\Theta_t u}^2\ \frac{dt}{t} 
	\leq \int_{0}^1 &\norm{\Theta_t(\iden - \Ppb_t) u}^2\ \frac{dt}{t} + 
	\int_{0}^1 \norm{(\Theta_t - \Pri_t \Av_t)\Ppb_t u}^2\ \frac{dt}{t} \\ &\qquad+ 
	\int_{0}^1 \norm{\Pri_t \Av_t(\Ppb_t - \iden)u}^2\ \frac{dt}{t} + 
	\int_{0}^1 \norm{\Pri_t \Av_t u}^2\ \frac{dt}{t} 
\end{aligned}
\end{equation}

The first three terms to the right of this expression can be 
handled relatively easily as the following lemma demonstrates. 

\begin{lemma}
We have that:
\begin{multline*}
\int_0^1 \norm{\Theta_t(\iden - \Ppb_t) u}^2\ \frac{dt}{t} + 
	\int_{0}^1 \norm{(\Theta_t - \Pri_t \Av_t)\Ppb_t u}^2\ \frac{dt}{t} \\ + 
	\int_{0}^1 \norm{\Pri_t \Av_t(\Ppb_t - \iden)u}^2\ \frac{dt}{t} 
	\lesssim \norm{A}_\infty^2 \norm{u}^2.
\end{multline*}
\end{lemma}
\begin{proof} 
For the first term, we estimate by noting that
$$\Theta_t(\iden - \Ppb_t) = \Theta_t t\Dirb \Qqb_t = (t \Pp_t \divv A_2 \Qqb_t)\Qqb_t,$$  
we obtain the required quadratic estimate 
using Proposition \ref{Prop:DivEst} to assert uniform
bounds for $t\Pp_t \divv A_2 \Qqb_t$ and by noting that  
$\Qqb_t$ satisfies quadratic estimates  \eqref{Eq:SA1}.  
The two remaining estimates are handled via
Propositions 5.4 and Proposition 5.10 in \cite{BMcR} 
with $S = \iden$. The versions of these 
propositions in our current context  can be obtained
exactly the way described in the proof of Proposition \ref{Prop:HarmAnal1}. 
\end{proof}

Thus, we have left with the last term in this expression, which 
we reduce to a Carleson measure estimate. That is, 
by Carleson's Theorem, the estimate of this term is
obtained by proving that 
$$ d\nu(x,t) = \modulus{\Pri_t(x)}^2\ \frac{d\mu(x)dt}{t}$$
is a Carleson measure. This is obtained 
if we prove for each cube $\Q \in \DyQ$, and for Carleson
regions $\CBox_{\Q} = \Q \times (0, \len(\Q))$,
\begin{equation}
\label{Eq:CarlEst}
\iint_{\CBox_{\Q}} \modulus{\Pri_t(x)}^2\ \frac{d\mu(x)dt}{t} \lesssim \norm{A}_\infty^2 \mu(\Q).
\end{equation} 

The estimate we perform here is more intricate
and involved than the Carleson measure estimate in Proposition 5.12 in \cite{BMcR}, 
and we provide
full details. First, observe the following important reduction.
\begin{lemma}
Suppose that  
for every cube $\Q \in \DyQ$ with $\len(\Q) \leq \met(\Q, \Sigma)$
the  Carleson estimate \eqref{Eq:CarlEst}  holds. Then, \eqref{Eq:CarlEst} holds for
every cube $\Q \in \DyQ$.
\end{lemma}
\begin{proof}
Fix  $\Q \in \DyQ^j$, with $j = \max\set{k_0, \jscale}$ (with $k_0$ coming from Lemma \ref{Lem:CanLem}), and
define the two sets  
\begin{align*}
\cA &= \set{\Q' \in \DyQ: \Q' \subset \Q\ \text{and}\  \met(\Q', \Sigma) \geq  \len(\Q')}, \\
\cB &= \set{\Q' \in \DyQ: \Q' \subset \Q\ \text{and}\  \met(\Q', \Sigma) < \len(\Q')}. \\
\end{align*} 
Now, consider the dyadic Whitney region  
$ \sW_{\Q'} = \Q' \times (\delta \len(\Q'), \len(\Q))$ 
so that 
$$ \CBox_{\Q} = \cbrac{ \Union_{\Q' \in \cA} \sW_{\Q'} } \union \cbrac{\Union_{\Q' \in \cB} \sW_{Q'}}.$$
Note that $\Q'' \subset \Q'$ and $\Q' \in \cA$ implies that $\Q'' \in \cA$.
Setting $\cA_{\mx}$ to be the maximal cubes in $\cA$, 
we obtain that
$$ \Union_{\Q' \in \cA} \sW_{\Q'} = \Union_{\Q' \in \cA_{\mx}} \CBox_{\Q'}.$$
On using the hypothesis, we obtain that
$$
\sum_{\Q' \in \cA_{\mx}} \iint_{\CBox_{\Q'}} \modulus{\Pri_t}^2\ \frac{d\mu dt}{t} 
	\lesssim \norm{A}_\infty^2 \sum_{\Q' \in \cA_{\mx}} \mu(\Q')
	\lesssim \norm{A}_\infty^2 \mu(\Q)$$
by the disjointedness of the cubes in $\cA_{\mx}$.

Next, note that from the off-diagonal decay of $\Theta_t$, 
we obtain that $\Theta_t: \Lp{\infty}(\cV) \to \Lp[loc]{2}(\cV)$, 
and reasoning as in \S5.2 in \cite{BMcR}, which comes from 
Corollary 5.3 in \cite{AKMC}, we have that
$$ \fint_{\Q'}\ \modulus{\Pri_t}^2 d\mu \lesssim \mu(\Q')$$
and therefore,  
$$ \int_{\sc_{\Q'}} \modulus{\Pri_t}^2\ \frac{d\mu dt}{t} 
	\lesssim \int_{\frac{\len(\Q')}{2}}^{\len(\Q')} \mu(\Q') \frac{dt}{t} 
	\lesssim \mu(\Q').$$
Now, fix $k > j$ and note that $\delta^k \leq \len(\Q)$ and 
for every cube $\Q' \in \cB_k = \cB \intersect \DyQ^k$, we
have that $\Q' \subset \set{x \in \Q: \met(x, \Sigma) \leq (C_1 + 1) \delta^{k}}$. 
 On invoking Lemma \ref{Lem:CanLem} with $s = \delta^{k} (C_1 + 1) \len(\Q)^{-1}$,
we obtain that 
$$ \mu(\Q') \lesssim \mu\set{x \in \Q: \met(x, \Sigma) \leq s \len(\Q)}
	\lesssim \frac{\delta^{k\tilde{\eta}}}{\len(\Q)^{\tilde{\eta}}} \mu(\Q) \lesssim \mu(\Q),$$
where the second inequality follows from $\delta^k \leq \len(\Q)$.   
 Note now that if $\Q' \in \cB$ and $\Q'' \subsetneq \Q'$
then $\len(\Q'') \leq \delta \len(\Q')$ and therefore,
$$\sW_{\Q'} = \Q' \times (\delta \len(\Q'), \len(\Q')) \intersect \Q'' \times (\delta \len(\Q''), \len(\Q'')) = \sW_{\Q''} = \emptyset,$$
and therefore 
$$ \sum_{\Q' \in \cB} \iint_{\sW_{\Q'}} \modulus{\Pri_t}^2\ \frac{d\mu dt}{t}
	\lesssim  \sum_{k > j} \sum_{\Q' \in \cB_k} \iint_{\sW_{\Q'}} \modulus{\Pri_t}^2\ \frac{d\mu dt}{t}  
	\lesssim \mu(\Q),$$
which completes the proof.
\end{proof}

We finally prove  \eqref{Eq:CarlEst} for the remaining
cubes $\Q$ bounded away from $\Sigma$.
\begin{proposition}
Suppose that $\met(\Q, \Sigma)  \geq \len(\Q)$.
Then, the Carleson measure estimate \eqref{Eq:CarlEst}
holds.
\end{proposition}
\begin{proof}
Fix $w \in \C^N$, let $f_{\Q}: \cM \to [0,1]$
with $\spt f_{\Q}$ compact, and $f_{\Q} = 1$ on $\Q$
 and $0$ outside $\Ball(x_{\Q}, 2 \len(\Q))$ with $\modulus{\conn f_{\Q}} \lesssim \len(\Q)^{-1}$. 
Define $w_{\Q}(x) = f_{\Q}(x) w^c(x) = f_{\Q}(x) w^i e_i(x)$
inside $B(x_{\ancester{\Q}}, \brad)$, the GBG trivialisation of $\Q$.
Note that, for $x \in \Q$  and $t \leq \scale$,  $\Av_tw_{\Q}(x) = w^c$.
Since the metric $\mh$ is uniformly comparable to the trivial 
metric inside this trivialisation, and using the 
facts we have just mentioned, 
$$ \iint_{\CBox_{\Q}} \modulus{\Pri_t}^2\ \frac{d\mu dt}{t} 
	\lesssim \sup_{\modulus{w}_\delta = 1} 
		\iint_{\CBox_{\Q}}\modulus{\Pri_t\Av_t w_{\Q}(x)}^2\ \frac{d\mu dt}{t}.$$
We split 
 \begin{multline*}
\iint_{\CBox_{\Q}}\modulus{\Pri_t\Av_t w_{\Q}(x)}^2\ \frac{d\mu dt}{t} \\ 
	\leq \iint_{\CBox_{\Q}}\modulus{(\Pri_t\Av_t - \Theta_t) w_{\Q}(x)}^2\ \frac{d\mu dt}{t}
	 	+  \iint_{\CBox_{\Q}}\modulus{\Theta_t w_{\Q}(x)}^2\ \frac{d\mu dt}{t}.
\end{multline*}
On following the exact same argument as in Proposition 5.11 in \cite{BMcR}, 
noting that this proof only requires that $\Theta_t$ satisfies
the off-diagonal estimates, we obtain 
that
$$\iint_{\CBox_{\Q}}\modulus{(\Pri_t\Av_t - \Theta_t) w_{\Q}(x)}^2\ \frac{d\mu dt}{t}
	\lesssim \norm{A}_\infty^2 \mu(\Q).$$

For the remaining part, let
$$ \Theta_t w_{\Q} = t \Pp_t \divv A_2 (\Ppb_t - \iden) w_{\Q}
	+ t \Pp_t \divv A_2 w_{\Q}.$$
We first obtain the required estimate on the second term.
For that, observe 
$w_{\Q} = 0$ near $\Sigma$ and hence,
$A_2 \omega_{\Q} \in \dom(\divv_{\mi})$.
 Using the identity $t\Pp_t \divv_{\mi} = (\Qq_t  +\imath t \Pp_t)
	\adj{(\conn(\imath\iden - \Dirp)^{-1})}$, 
we estimate 
$$ \int_{0}^{\len(Q)} \norm{t\Pp_t \divv A_2 w_{\Q}}^2\ \frac{dt}{t}
	\lesssim \norm{ \adj{(\conn(\imath\iden - \Dirp)^{-1})} A_2 w_{\Q}}^2 \lesssim \norm{A}_\infty^2 \mu(\Q).$$

To estimate the remaining term, 
we note that   $t\Pp_t \divv A_2 (\Ppb_t -\iden) w_{\Q} = - t \Pp_t \divv A_2 \Qqb_t (t\Dirb w_{\Q})$  
and so   by Proposition \ref{Prop:DivEst}  
$$ \norm{t\Pp_t \divv A_2 (\Pp_t -\iden) w_{\Q}}^2
 	\lesssim  t^2 \norm{A}_\infty^2 \norm{\Dirb w_{\Q}}^2
	\lesssim  t^2 \norm{A}_\infty^2 \norm{\conn w_{\Q}}^2  
	\lesssim t^2 \norm{A}_\infty^2 \frac{1}{\len(\Q)^2} \mu(\Q).$$
Therefore,
\begin{multline*}
  \iint_{\CBox_{\Q}} \modulus{t\Pp_t \divv A_2 (\Pp_t -\iden) w_{\Q}}^2\ \frac{d\mu dt}{t} 
	\leq \int_{0}^{\len(\Q)}\norm{t\Pp_t \divv A_2 (\Ppb_t -\iden) w_{\Q}}^2 \frac{dt}{t}  \\
	\lesssim \norm{A}_\infty^2 \mu(\Q)  \int_{0}^{\len(\Q)} \frac{t}{\len(\Q)^2} dt 
	\lesssim \norm{A}_\infty^2 \mu(\Q),
\end{multline*}
which establishes the conclusion.
\end{proof}

\begin{proof}[Proof of Theorem \ref{Thm:Main}]
On combining the estimates in \S\ref{Sec:Harm2}
and  Proposition  \ref{Prop:HarmAnal1},  the hypothesis
of Proposition \ref{Prop:FinalRed} is satisfied. 
This proves Theorem 2.1.
\end{proof}

\bibliographystyle{amsplain}
\def\cprime{$'$}
\providecommand{\bysame}{\leavevmode\hbox to3em{\hrulefill}\thinspace}
\providecommand{\MR}{\relax\ifhmode\unskip\space\fi MR }
\providecommand{\MRhref}[2]{%
  \href{http://www.ams.org/mathscinet-getitem?mr=#1}{#2}
}
\providecommand{\href}[2]{#2}

\setlength{\parskip}{0mm}

\end{document}